\documentclass[a4paper,12pt,oneside]{article}
\usepackage[margin=1in]{geometry}
\pdfoutput=1

\RequirePackage{amsmath,amsthm, amssymb,mathtools}
\RequirePackage{natbib}
\RequirePackage{algorithm2e}
\RequirePackage{tikz}
\RequirePackage{pgfplots}
\RequirePackage{float}
\RequirePackage{enumitem}

\usepackage{libertine}
\usepackage{courier}
\usepackage{mathpazo}
\usepackage{xcolor}
\usepackage{pifont}
\usepackage{authblk}
\usepackage[backref=page]{hyperref}
\RequirePackage{nicematrix}
\NiceMatrixOptions{
code-for-first-row = \tiny ,
code-for-last-row = \tiny ,
code-for-first-col = \tiny ,
code-for-last-col = \tiny
}

\usepackage{titlesec}

\definecolor{darkbrown}{RGB}{170,100,0} 
\definecolor{darkdarkbrown}{RGB}{110,70,0} 
\titleformat{\section}
{\color{darkdarkbrown}\normalfont\Large\bfseries}
{\color{darkdarkbrown}\thesection}{1em}{}
\titleformat{\subsection}
{\color{darkdarkbrown}\normalfont\Large\bfseries}
{\color{darkdarkbrown}\thesubsection}{1em}{}

\hypersetup{
	colorlinks=true,
  linkcolor=darkbrown,
  citecolor=darkbrown,
  filecolor=darkbrown,
  urlcolor=darkbrown
}

 \renewcommand*{\backrefalt}[4]{%
    \ifcase #1%
     \or (page:~#2)%
     \else (pages:~#2)%
    \fi%
    }

    \makeatletter
\def\@fnsymbol#1{\ensuremath{\ifcase#1\or \dagger\or \ddagger\or
   \mathsection\or \mathparagraph\or \|\or **\or \dagger\dagger
   \or \ddagger\ddagger \else\@ctrerr\fi}}
    \makeatother

\newtheorem{theorem}{Theorem}[section]
\newtheorem{lemma}{Lemma}[section]

\newtheorem{proposition}{Proposition}[section]
\newtheorem{example}{Example}[section]
\newtheorem{corollary}{Corollary}[section]
\newtheorem{remark}{Remark}[section]

\providecommand{\keywords}[1]
{
  \small	
  \textbf{\textit{Keywords---}} #1
}
\newcommand{\pred}[1]{\boldsymbol{1}\left\{#1\right\}}
\newcommand{\N}{\mathbb{N}}

\newcommand{\trn}{^\intercal}
\newcommand{\inv}{^{-1}}

\DeclareMathOperator*{\argmin}{arg\,min}
\DeclareMathOperator*{\argmax}{arg\,max}
\newcommand{\abs}[1]{\left| #1 \right|}
\newcommand{\nrm}[1]{\left\Vert #1 \right\Vert}
\newcommand{\tv}[1]{\nrm{#1}_{\mathsf{TV}}}

\newcommand{\esttmix}{\widehat{t}_{\mathsf{mix}}}

\newcommand{\sg}{\gamma}
\newcommand{\bbE}{\mathbb{E}}

\newcommand{\pssg}{\sg_{\mathsf{ps}}}

\newcommand{\tmix}{t_{\mathsf{mix}}}
\newcommand{\trel}{t_{\mathsf{rel}}}
\newcommand{\PR}[2][]{\mathbb{P}_{#1}\left( #2 \right)}
\newcommand{\E}[2][]{\mathbb{E}_{#1}\left[ #2 \right]}

\newcommand{\as}{\stackrel{a.s.}{\to}}
\newcommand{\eps}{\varepsilon}
\DeclareMathOperator{\Bin}{Bin}

\DeclareMathOperator{\Uniform}{Uniform}
\DeclareMathOperator{\Dir}{Dir}

\newcommand{\estpssgS}{\widehat{\sg}_{\mathsf{ps} [S]}}

\newcommand{\Nmin}{N_{\min}}
\renewcommand{\ln}{\log}

\newcommand{\eqdef}{\doteq}
\newcommand{\pimin}{\pi_\star}

\newcommand{\rev}{^{\star}}
\newcommand{\bigO}{\mathcal{O}}
\newcommand{\hide}[1]{}

\newcommand{\genstar}{\star}
\newcommand{\cc}{\kappa}

\newcommand{\tkgen}{{s_{\star}}}

\newcommand{\tgencc}{\cc_{\star}}

\newcommand{\estcc}{\widehat{\cc}}

\newcommand{\esttgencc}{\widehat{\cc}_{\genstar[S]}}
\newcommand{\esttgenccadapt}{\widehat{\cc}_{\genstar}}
\newcommand{\esttgenccexpl}{\widehat{\cc}_{\genstar\ceil{2/\eps}}}

\newcommand{\esttgenccexplest}{\widehat{\cc}_{\genstar[\widehat{S}]}}
\newcommand{\esttgenccexplestbb}{\widehat{\cc}_{\genstar[S]}}
\newcommand{\esttgenccpea}{\widehat{\cc}_{\genstar}^{+}}

\newcommand{\tgenccS}{\cc_{\genstar[S]}}

\newcommand{\calX}{\mathcal{X}}

\newcommand{\calP}{\mathcal{P}}
\newcommand{\calW}{\mathcal{W}}
\newcommand{\calD}{\mathcal{D}}
\newcommand{\R}{\mathbb{R}}

\DeclarePairedDelimiter\ceil{\lceil}{\rceil}
\DeclarePairedDelimiter\floor{\lfloor}{\rfloor}
\newcommand{\set}[1]{\left\{ #1 \right\}}

\numberwithin{equation}{section}

\title{Empirical and Instance--Dependent \\ Estimation of Markov Chain and \\ Mixing Time}

\author{Geoffrey Wolfer \thanks{email: geoffrey.wolfer@riken.jp \\ The author is supported by the Special Postdoctoral Researcher Program (SPDR) of RIKEN.} \\ Center for AI Project, RIKEN}

\date{\today}

\begin{document}

\maketitle

\begin{abstract}
We address the problem of estimating the mixing time of a Markov chain from a single trajectory of observations. Unlike most previous works which employed Hilbert space methods to estimate spectral gaps, we opt for an approach based on contraction with respect to total variation. Specifically, we estimate the contraction coefficient introduced in \citet{pmlr-v117-wolfer20a}, inspired from Dobrushin's. This quantity, unlike the spectral gap, controls the mixing time up to strong universal constants and remains applicable to non--reversible chains. We improve existing fully data--dependent confidence intervals around this contraction coefficient, which are both easier to compute and thinner than spectral counterparts. Furthermore, we introduce a novel analysis beyond the worst--case scenario by leveraging additional information about the transition matrix. This allows us to derive instance--dependent rates for estimating the matrix with respect to the induced uniform norm, and some of its mixing properties.

\end{abstract}


\clearpage
\tableofcontents
\clearpage

\section{Introduction}
\label{section:introduction}
This study delves into the intricacies of constructing non--trivial high confidence intervals for the transition matrix and the mixing time of a finite state time--homogeneous ergodic Markov chain. Challenges arise when only a single long trajectory of states, $X^m = X_1, X_2, \dots, X_m$, is available for observation, without access to a restart mechanism.
The problem is motivated by the following applications.
\begin{enumerate}
    \item Probably Approximately Correct (PAC)--type learning problems, which assume data sampled from a Markovian process,
and where generalization guarantees frequently depend on the \emph{a priori} unknown mixing properties of the chain \citep{kuznetsov2017generalization}. 
\item Markov Chain Monte Carlo (MCMC) diagnostics for non--reversible Markov chains, which may exhibit superior mixing
properties or asymptotic variance compared to their reversible counterparts \citep{bierkens2016non}.
\item Reinforcement learning, where bounds on the mixing
time of Markov decision processes \citep{ortner2020regret} are routinely assumed.
\end{enumerate}
We encourage the reader to refer to the related work sections of \cite{hsu2019, wolfer2023improved} for a more comprehensive collection of references regarding the aforementioned problems. In this paper, we are particularly interested in obtaining confidence intervals and procedures which could benefit from additional structure, or known properties of the transition matrix, e.g. sparsity of its connection graph. By doing do, our goal is to provide a beyond worst--case analysis for two central problems in statistics.

\subsection{Main results and outline}

We summarize our main technical contributions as follows.
\begin{enumerate}
    \item We construct a fully--empirical confidence interval for the task of estimating the transition matrix of a Markov chain from a single trajectory of $m$ observations (Theorem~\ref{theorem:fully-empirical-transition-kernel-learning}). The interval decays at a rate $\widetilde \bigO(\sqrt{m})$ and can leverage sparsity properties of the transition matrix, providing an empirical counterpart of the known minimax results of \citet{wolfer2021}.
    Crucially, no additional assumptions have to be made on the unknown Markov chain; the interval will be naturally thinner in favorable cases.
    \item Improving upon \citet[Theorem~2]{pmlr-v117-wolfer20a}, we obtain fully--empirical confidence intervals for the mixing time, which depend more delicately on the data. In contrast with prior work, where the estimator has a hard--coded choice of an approximation parameter, the proposed estimator (Theorem~\ref{theorem:empirical-estimation-tmix}) is completely determined from the data.
    \item Devising an amplification method, we show that the minimax trajectory length for the problem of estimating the contraction coefficient $\kappa_\star$ of \citet[(9)]{pmlr-v117-wolfer20a} to constant multiplicative error is upper bound by
    \begin{equation*}
        \widetilde{\bigO} \left( \Xi(P) + \frac{1}{(1 - \kappa_\star) \pi_\star} \right),
    \end{equation*}
    where $\pi_\star$ is the minimum stationary probability \eqref{eq:minimum-stationary-probability} of $P$, and $\Xi(P)$ is a quantity that depends delicately on $P$ (Theorem~\ref{theorem:learn-tgencc-ub-relative}). This improves the previous upper bound of \citet[Theorem~6]{pmlr-v117-wolfer20a} by shaving off at least a factor $1/(1 - \kappa_\star)$.
    \item Along the way, we provide a tightness analysis for \citet{pmlr-v117-wolfer20a}, derive additional properties of $\kappa_\star$ as well as a minimax upper bound for the problem of estimating the transition matrix of a skipped Markov chain (Theorem~\ref{theorem:learning-s-step-transition-matrix}).
\end{enumerate}

\paragraph{Outline}

In Section~\ref{section:generalized-contraction-coefficient}, we begin by recalling the definition of the contraction coefficient $\kappa_\star$ \eqref{eq:generalized-contraction-coefficient} of a transition matrix $P$ proposed in \citet[(9)]{pmlr-v117-wolfer20a}, and the main property that $1/(1 - \kappa_\star)$ controls the mixing time of $P$ up to tight universal constants, unlike spectral quantities that lead to a logarithmic gap \eqref{eq:mixing-time-vs-relaxation-time}. We then turn our attention to examining additional properties of $\kappa_\star$, e.g. in terms of optimality (Theorem~\ref{theorem:control-tau-with-kgen-and-xi-converse}), or in relation to skipped chains (Property~\ref{proposition:kgen-of-power}).
In Section~\ref{section:empirical-estimation-matrix}, we tackle the problem of estimating $P$ with respect to the $\ell_\infty$ matrix norm. We construct an estimator $\widehat{P}$ together with a confidence interval that is computable from the observed data and traps the true $P$ with high--probability. We show that this approach is consistent with known minimax rates, supersedes currently known empirical bounds, is able to leverage fine--grained properties of the transition matrix to speed--up convergence guarantees, and we provide a few examples.
In Section~\ref{section:empirical-estimation-tmix}, we improve the  estimator of \citet{pmlr-v117-wolfer20a} for $\kappa_\star$ based on the given trajectory observed at different skipping rates.
We show how to choose to the maximum skipping rate parsimoniously and adaptively from the data, and perform an asymptotic analysis of the interval.
As a direct consequence, we construct an estimator and fully empirical confidence intervals for the mixing time of $P$ itself.
In Section~\ref{section:minimax-estimation}, we analyze the sample complexity of estimating the multi--step transition kernel $P^s$, the coefficient $\kappa_\star$ and the mixing time through the lens of minimax theory. We show that the properties of $\kappa_\star$, despite being weaker than that of the absolute spectral gap in the reversible setting, still enable us to perform amplification of the estimator and convert estimation rates from additive to multiplicative accuracy.
Finally, in Section~\ref{section:family-special}, we inspect the properties of a family of three--state Markov chains, and show that the skipping rate that achieves the coefficient $\kappa_\star$ can be arbitrarily large.

\subsection{Notation and background}
\label{section:definition}

The set $\N$ refers to the natural integers $\set{1, 2, \dots}$
and for $n \in \N$, we write $[n] = \set{1, 2, \dots, n}$.
Unless a base is specified, 
$\log$ denotes the natural logarithm.
We henceforth fix a finite set $\calX$, 
and define $\calP(\calX)$, the simplex of all distributions ---seen as row vectors--- over $\calX$. 
For $x \in \calX$, we denote by $e_x$ the vector such that $\forall x' \in \calX, e_x(x') = \pred{x = x'}$.
For $v \in \R^\calX$, and $p \in \R_+$, we write $\nrm{v}_p \eqdef \left(\sum_{x \in \calX} \abs{v(x)}^p\right)^{1/p}$,
which corresponds to the $\ell_p$ norm when $p \geq 1$.
For $\mu \in \calP(\calX)$ and $n \in \N$, we also consider the $n$--fold product $\mu^{\otimes n}$, i.e.
$X_1, \dots, X_n \sim \mu^{\otimes n}$ is a shorthand for $(X_t)_{t \in [n]}$ being all mutually independent and such that
$\forall t \in [n], X_t \sim \mu$.
For $(\mu, \nu) \in \calP(\calX)^2$,
we define the total variation distance (TV) between $\mu$ and $\nu$ in terms of the $\ell_1$ norm as
\begin{equation}
\label{eq:total-variation}
\tv{\mu -\nu} \eqdef \frac{1}{2} \nrm{\mu -\nu}_1.
\end{equation}
We consider time--homogeneous Markov chains
$$X^{\infty} = X_1, X_2, \dots, X_t, \dots \; \sim (\mu, P),$$ 
with initial distribution $\mu \in \calP(\calX)$,
and row--stochastic transition matrix $P \colon \calX \times \calX \to [0,1]$, $\forall x \in \calX, \nrm{e_xP}_1 = 1$.
The set of all row--stochastic matrices over $\calX$ is denoted by $\calW(\calX)$.
A Markov chain is called ergodic when its transition matrix $P$ is primitive, 
i.e. $\exists p \in \N, P^p > 0$ entry--wise. 
In this case, $P$ has a unique stationary distribution $\pi$ verifying $\pi P = \pi$,
and the minimum stationary probability,
\begin{equation}
\label{eq:minimum-stationary-probability}
\pimin \eqdef \min_{x \in \calX} \pi(x),
\end{equation}
satisfies $\pimin > 0$. 
Measuring the distance to stationarity in total variation,
\begin{equation}
\label{eq:distance-to-stationarity}
h(t) \eqdef \sup_{\mu \in \calP(\calX)} \tv{\mu P^t - \pi},
\end{equation}
the chain converges to $\pi$ in the sense where $\lim_{t \to \infty} h(t) = 0$.
More quantitatively, for $\xi \in (0, 1/2)$, the mixing time of $P$ is defined by
\begin{equation}
\label{eq:mixing-time}
\tmix(\xi) \eqdef \argmin_{t \in \N} \set{h(t) < \xi},
\end{equation}
and by convention $\tmix \eqdef \tmix(1/4)$. 
We refer the reader to \cite[Chapter~4]{levin2009markov} 
for a comprehensive introduction to Markov chain mixing.
The definition of elements specific to contraction methods is deferred to
Section~\ref{section:generalized-contraction-coefficient} for clarity of the exposition.

\subsection{Related work}
The mixing time estimation problem from a single trajectory of observations has been investigated, using spectral methods, by \citet{hsu2015mixing, levin2016estimating, hsu2019, combes2019computationally} in the time--reversible setting, and by \citet{pmlr-v99-wolfer19a, wolfer2023improved} in the more prevalent non--reversible case.
All methods in this body of work hinge around the fact that one can effectively upper and lower bound the mixing time by another quantity, termed the (pseudo--)relaxation time $\trel$ of the transition matrix,
\begin{equation}
\label{eq:mixing-time-vs-relaxation-time}
c_1 \trel \leq \tmix \leq c_2 \trel \log \frac{1}{\pimin},
\end{equation}
where $c_1$ and $c_2$ are two universal constants.
When $P$ is reversible, the relaxation time is simply defined as the inverse of the absolute spectral gap $\gamma_\star$, which is the difference between the two eigenvalues of largest magnitude. The curious reader can find a proof of \eqref{eq:mixing-time-vs-relaxation-time} in \citet[Theorem~12.4, Theorem~12.5]{levin2009markov}.
The absence of reversibility raises serious challenges for spectral methods. Specifically, for $\trel = 1/\gamma_\star$, the logarithmic term $\log 1 / \pi_\star$ in the upper bound in \eqref{eq:mixing-time-vs-relaxation-time} must be replaced with a linear dependency in the state space size $\abs{\calX}$ \citep[Theorem~1.2]{jerison2013general}.
In this setting, the pseudo--relaxation time \citep{kamath2016estimation} can be better defined as the inverse of the pseudo--spectral gap $\pssg$, introduced by \citet{paulin2015concentration}
and both inequalities in \eqref{eq:mixing-time-vs-relaxation-time} hold with $\trel = 1 / \pssg$ \citep[Proposition~3.4]{paulin2015concentration}.
It is instructive to briefly inspect the definition of $\pssg$,
\begin{equation}
    \label{eq:pseudo-spectral-gap}
    \pssg \eqdef \max_{k \in \N} \set{\frac{1}{k}\sg_\star\left((P \rev)^kP^k\right)}.
\end{equation}
$P \rev$ is the time--reversal of $P$, which governs the dynamics of the stationary Markov chain observed in reverse time. The product $P^\star P$, termed multiplicative reversiblization \citep{fill1991eigenvalue} of $P$, is by construction reversible and its absolute spectral gap can be used to bound the mixing time of $P$, albeit pessimistically.
For a more detailed account of spectral methods, We refer the reader to the related work section of \citet{wolfer2023improved}.

One severe limitation of spectral approaches, which is also raised in \citet[p.511--512]{rabinovich2020function}, is the existence of an unclosable \citep{jerison2013general} logarithmic gap $\log 1 / \pi_\star$ 
in \eqref{eq:mixing-time-vs-relaxation-time}.
The direct consequence is that knowledge of $\trel$ even down to arbitrary small error, does not translate into a reliable estimate for $\tmix$ as $\abs{\calX} \rightarrow \infty$. 
In fact, it is currently believed that the 
problem of estimating $\tmix$ to constant multiplicative error is more challenging than $\trel$ (see e.g. concluding remarks of \cite[Conclusion]{combes2019computationally}).

Inspired by Dobrushin and Paulin, \citet{pmlr-v117-wolfer20a} departed from spectral methods, and proposed a tighter proxy, up to universal constants, based on higher--order contraction coefficients ---reported in \eqref{eq:generalized-contraction-coefficient}--- as well as an algorithm to estimate it from a single trajectory.
The approach therein compares favorably with spectral methods in that it treats all ergodic chains, regardless of reversibility.

\section{The contraction coefficient \texorpdfstring{$\tgencc$}{}}
\label{section:generalized-contraction-coefficient}
The Dobrushin contraction coefficient, also known as Dobrushin ergodic coefficient 
\citep{dobrushin1956central}, \citep[Definition~7.1]{bremaud99} of a Markov chain with transition matrix $P$ is defined as
\begin{equation}
\label{eq:dobrushin-coefficient}
\cc \eqdef \max_{\mu, \nu} \tv{ \mu P - \nu P},
\end{equation}
where the maximum is taken over pairs of distributions on $\calX$, and the 
term contraction refers to the property \citep[Corollary~7.1]{bremaud99} 
that two probability measures can only be brought closer together by the action of the kernel; more formally,
\begin{equation}
\label{eq:contraction-property}
\forall (\mu,\nu) \in \calP(\calX)^2, \tv{(\mu-\nu)P } \le \cc \tv{\mu-\nu}.
\end{equation}
Contraction in the sense of Dobrushin is a special case of coarse Ricci curvature \citep{ollivier2009ricci},
where the metric taken on $\calX$ is the discrete metric, and the Wasserstein 
distance between distributions reduces to total variation. In the case where $\cc < 1$, 
we can readily upper bound the mixing time as follows,
\begin{equation}
\label{eq:bubley-dyer-path-coupling-bound}
\begin{split}
\tmix(\xi) \leq \frac{\ln{\xi}}{\ln{(1 - \cc)}}.
\end{split}
\end{equation}
However, an ergodic Markov chain could still be non--contractive, i.e. $\cc = 1$. 
In this case, the aforementioned method fails to yield convergence rates.
Furthermore, even for a contractive chain, it is not possible to derive lower 
bounds on the mixing time in terms of $\kappa$, and the upper bound 
at \eqref{eq:bubley-dyer-path-coupling-bound} may be pessimistic.
It is natural to consider multi--step chains, where for $s \in \N$, 
\begin{equation*}
\begin{split}
X_{[s]}^{\infty} \eqdef X_{1}, X_{1 + s}, X_{1 + 2s}, \dots, X_{1 + ts}, \dots
\end{split}
\end{equation*}
and define (see e.g. \citet{paulin2016mixing}) the contraction coefficient of the chain with skipping rate $s$ to be
\begin{equation}
\label{eq:dobrushin-coefficient-skipped}
\cc_s \eqdef \max_{\mu, \nu} \tv{\mu P^s - \nu P^s}.
\end{equation}
For $s \in [m-1]$, we define
\begin{equation*}
\begin{split}
X^m_{[s]} \eqdef X_{1}, X_{1 + s}, X_{1 + 2s}, \dots, X_{1 + \floor{(m-1)/s}s} \sim (\mu, P^s),
\end{split}
\end{equation*}
for the $s$--skipped chain observed from an original trajectory of length $m$.
For $\xi \in (0, 1/2)$, we write
\begin{equation}
\label{eq:tau-definition}
\tau(\xi) \eqdef \argmin_{s \in \N} \set{\cc_s < \xi}, \qquad \tau \eqdef \tau(1/4),
\end{equation}
which are closely related to $\tmix(\xi)$ and $\tmix$ \citep[Lemma~4.10]{levin2009markov}, as by the triangle inequality,
\begin{equation}
\label{eq:link-tmix-tau}
\tau(2 \xi) \leq \tmix(\xi) \leq \tau(\xi).
\end{equation}
\citet[(9)]{pmlr-v117-wolfer20a} proposed the generalized contraction 
coefficient $\tgencc$ of the ergodic chain $P$:
\begin{equation}
\label{eq:generalized-contraction-coefficient}
\tgencc \eqdef 1 - \max_{s \in \N} \set{\frac{1 - \cc_s}{s}},
\end{equation}
constructed in a similar spirit to \citeauthor{paulin2015concentration}'s pseudo--spectral gap of an ergodic Markov chain, based on a collection of 
higher--order multiplicative reversiblizations \citep{fill1991eigenvalue}.
It will be convenient to write $\tkgen$ the smallest integer 
\footnote{The existence of $\tkgen$ is guaranteed by the observation that 
$s \mapsto \frac{1 - \cc_s}{s} \in [0, 1/s]$.} such that 
$\tgencc = 1 - \frac{1 - \cc_{\tkgen}}{\tkgen}.$
We immediately observe that even for a non--contractive ergodic chain, 
it always holds that $\tgencc < 1$.
Indeed, by definition of ergodicity (Section~\ref{section:definition}), 
there exists $p \in \N$ and $\alpha > 0$ such that
$P^p \geq \alpha$ (entry--wise), thus $\cc_p \leq 1 - \abs{\calX} \alpha$,
and $\tgencc \leq 1 - \frac{\abs{\calX} \alpha}{p} < 1$.
We first recall, in Theorem~\ref{theorem:control-tau-with-kgen-and-xi-alt2020} and Corollary~\ref{corollary:control-tmix-with-kgen-alt2020}, the stronger result that $\frac{1}{1 - \tgencc}$ effectively traps 
$\tau$ and $\tmix$ up 
to multiplicative universal constants.

\begin{theorem}[{\citealt[Theorem~1]{pmlr-v117-wolfer20a}}]
\label{theorem:control-tau-with-kgen-and-xi-alt2020}
For any $\xi \in (0,1/2)$, and for any ergodic transition matrix $P$ over a space $\calX$,
	there exists a universal constant $c \leq 1$ such that
\begin{equation*}
\begin{split}
\frac{1 - \xi}{1 - \tgencc} \leq \tau(\xi) \leq \frac{\ln (c e / \xi)}{1 - \tgencc},
\end{split}
\end{equation*}
where $\tgencc$ is defined at \eqref{eq:generalized-contraction-coefficient}, 
$\tau(\xi)$ is defined at \eqref{eq:tau-definition},
and $e = 2.71828\dots$ is the Euler number.
\end{theorem}

\begin{corollary}[{\citealt[Theorem~1]{pmlr-v117-wolfer20a}}]
\label{corollary:control-tmix-with-kgen-alt2020}
Let $P$ be any ergodic transition matrix over a space $\calX$ with mixing time $\tmix$.
Then there exist two universal constants 
$\underline{c} \geq 1/2$ and $ \overline{c} \leq \ln{4e} = 2.38629\dots$ such that
\begin{equation*}
\begin{split}
\frac{\underline{c}}{1 - \tgencc} \leq \tmix \leq \frac{\overline{c}}{1 - \tgencc}.
\end{split}
\end{equation*}
\end{corollary}
We complement the above result by providing converse statements, demonstrating the tightness of Theorem~\ref{theorem:control-tau-with-kgen-and-xi-alt2020}.
\begin{theorem}
\label{theorem:control-tau-with-kgen-and-xi-converse}
Let $\xi \in (0,1/2)$.
\begin{enumerate}
\item[$(a)$] The universal constant $c$ in Theorem~\ref{theorem:control-tau-with-kgen-and-xi-alt2020} satisfies $c \geq 1/2$.
\item[$(b)$] There exists an ergodic transition matrix $P$ 
such that $\tau(\xi) = \frac{1 - \xi}{1 - \tgencc}$.
\item[$(c)$] There exists an ergodic transition matrix $P$ 
such that $\tau(\xi) \geq \frac{\ln (1/\xi)}{1 - \tgencc}$.
\end{enumerate}
\end{theorem}
In Section~\ref{section:numerical-tightness-bound}, we illustrate numerically the tightness of Theorem~\ref{theorem:control-tau-with-kgen-and-xi-alt2020} and Theorem~\ref{theorem:control-tau-with-kgen-and-xi-converse}.
It is important to note that the inequalities of Theorem~\ref{theorem:control-tau-with-kgen-and-xi-alt2020} and 
Corollary~\ref{corollary:control-tmix-with-kgen-alt2020} do not involve $\abs{\calX}$ or $\pimin$, thus
contrary to the bound at \eqref{eq:mixing-time-vs-relaxation-time} obtained with spectral methods, 
they do not degrade as $\abs{\calX} \rightarrow \infty$.
The following proposition establishes connections between the generalized contraction
coefficients of a chain and its multi--step version.\\

\begin{proposition}
\label{proposition:kgen-of-power}
Let $p \in \N$.
Let $P$ be an ergodic transition matrix over a space $\calX$,
and $P^p$ its $p$th power.
Let $\tmix$ and $\tgencc$ (resp. $\tgencc^{(p)}$)
pertain to $P$ (resp. $P^p$). 
\begin{enumerate}[label=$(\roman*)$]
    \item It holds that
    $$1 - \tgencc^{(p)} \leq p (1 - \tgencc).$$
    \item When $p \geq \tmix$, 
    $$1 - \tgencc^{(p)} \geq 1/2.$$
    \item There exists a universal constant $c_\star = \frac{1}{4 \ln {4e}} = 0.10476\dots$ such that when $p \leq \tmix$, 
    $$1 - \tgencc^{(p)} \geq c_\star p (1 - \tgencc).$$
\end{enumerate}
\end{proposition}

Recall that for a reversible $P$, the pseudo--spectral gap verifies, 
$\pssg = \sg(P \rev P)$, i.e. the maximum is achieved for $k = 1$ in \eqref{eq:pseudo-spectral-gap} \citep[Lemma~15]{pmlr-v99-wolfer19a}.
We now emphasize that this fact does not translate to $\tgencc$ with respect to the contractivity property. 
In fact, even for a contractive chain over only three states (see the family constructed in Section~\ref{section:family-special}), the skipping rate $s_\star$ that achieves the maximum in the definition of $\tgencc$ can be arbitrarily large.
\begin{theorem}
\label{theorem:arbitrary-large-s}
For any $s_0 \in \N$, there exists an ergodic transition matrix $P$, such that 
$$\argmax_{s \in \N} \set{ \frac{1 - \cc_s}{s}} \subset [s_0, \infty).$$
\end{theorem}
\begin{proof}
See Lemma~\ref{lemma:special-category-chain-requires-arbitray-large-s}.
\end{proof}

\section{Empirical estimation of the transition matrix}
\label{section:empirical-estimation-matrix}
As an introductory result, we first derive an empirical confidence interval
for the problem of estimating the transition matrix with respect to the (matrix) uniform norm.
Namely, for a confidence parameter $\delta \in (0, 1)$, and a single trajectory of observations $X^m = X_1, \dots, X_m$, 
we construct an estimator $\widehat{P}(X^m)$ for $P$ and an interval width $\widehat{W}_{\delta}(X^m) \subset [0, 2]$, 
that verifies
\begin{equation*}
\begin{split}
\PR{\nrm{ \widehat{P}(X^m) - P }_\infty \leq \widehat{W}_{\delta}(X^m)} &\geq 1 - \delta,
\end{split}
\end{equation*}
and such that $\widehat{W}_{\delta}$ is non--trivial, i.e. $\abs{\widehat{W}_{\delta}(X^m)} \as 0$.
It is known \citep[Theorem~3.1]{wolfer2021}, that for a trajectory length of at least 
\begin{equation*}
    m \geq \frac{c}{\pimin} \max \set{\frac{1}{\eps^2 } \max \set{ \abs{\calX}, \log \frac{1}{\eps \delta} }, \frac{1}{\pssg} \log \frac{\abs{\calX} \nrm{\mu / \pi}_{\pi}}{\delta}},
\end{equation*}
where 
\begin{equation}
\label{eq:pi-norm}
\nrm{\mu/\pi}_{\pi} \eqdef \sum_{x \in \calX} \mu(x)^2 /\pi(x) \leq 1/\pimin, 
\end{equation}
$\pi, \pi_\star, \pssg$ pertain to the unknown transition matrix $P$, $\eps$ is the precision parameter, $\mu$ is the initial distribution, and $c$ is a universal constant, it holds that
\begin{equation*}
\begin{split}
\PR{\nrm{ \widehat{P} - P }_\infty \leq \eps} &\geq 1 - \delta.
\end{split}
\end{equation*}
Furthermore, the above sample complexity is known to be minimax optimal (up to logarithmic factors), and achieved for some classes of reversible Markov chains \citep[Theorem~3.2]{wolfer2021}.
However, while irreducibility makes it possible to bound $\nrm{\mu / \pi}_{\pi}$ in terms of $1 / \pi_\star$, the necessary trajectory length still depends on the unknown parameters $\pi_\star$ and ---perhaps ironically--- the mixing parameter $\pssg$.
 We show that we can bypass the requirement of finding estimates for these two quantities for this problem, and that the natural estimator based on a smoothed tally matrix enjoys fully empirical confidence intervals.
We define the natural counting random variables,
\begin{equation}
\label{eq:counting-random-variables-no-skip}
\begin{split}
N_{x} &\eqdef \sum_{t=1}^{\floor{m-1}} \pred{X_{t} = x}, \;\; \Nmin \eqdef \min_{x \in \calX} N_x, \;\; N_{xx'} \eqdef \sum_{t=1}^{\floor{m-1}} \pred{X_{t} = x, X_{1 + t} = x'}.
\end{split}
\end{equation}
Given a fixed smoothing vector $\lambda \in (0, \infty)^{\calX}$, 
we then construct a smoothed estimator for $P$,
\begin{equation}
\label{eq:empirical-transition-matrix-no-skip}
\begin{split}
\widehat{P}_{\lambda} &\eqdef \sum_{x,x' \in \calX} \frac{N_{xx'} + \lambda(x) }{N_x + \lambda(x) \abs{\calX}} \; e_x \trn e_{x'}.
\end{split}
\end{equation}
Notice that the entries of $\widehat{P}_{\lambda}$ are all positive,
hence the estimator is a almost surely contractive.
Moreover, the stationary distribution of $\widehat{P}_{\lambda}$ is given by
\begin{equation*}
\begin{split}
\widehat{\pi}_{\lambda} \eqdef \sum_{x \in \calX} \frac{N_x + \lambda(x) \abs{\calX}}{m + \abs{\calX} \nrm{\lambda}_1 } e_x.
\end{split}
\end{equation*}
The below--stated theorem provides a fully empirical confidence interval for 
estimating the transition kernel of a Markov chain with respect 
to the $\ell_\infty$ operator norm.
\begin{theorem}
\label{theorem:fully-empirical-transition-kernel-learning}
Let $P$ be an ergodic transition matrix over a finite space $\calX$, 
$\mu \in \calP(\calX)$, $\delta \in (0,1)$, 
$\lambda \in (0,\infty)^{\calX}$ a smoothing vector, 
$m \in \N$, and $X^m = X_1, X_2, \dots, X_m \sim (\mu, P)$.
With probability at least $1 - \delta$, it holds that
\begin{equation*}
\begin{split}
\nrm{\widehat{P}_{\lambda}(X^m) - P}_\infty &\leq \widehat{W}_{\lambda, \delta}(X^m),\\
\end{split}
\end{equation*}
with
\begin{equation}
\label{equation:w-confidence-interval}
\begin{split}
\widehat{W}_{\lambda, \delta}(X^m) &\eqdef 2 \max_{x \in \calX} \set{ \frac{  \sum_{x' \in \calX} \sqrt{N_{xx'}} + (3/\sqrt{2}) \sqrt{N_x} \sqrt{\ln {(2 m \abs{\calX} /\delta)}}  + \lambda(x)\abs{\calX} }{N_x + \lambda(x)\abs{\calX} }},
\end{split}
\end{equation}
and where $N_x$ and $N_{xx'}$ are defined at \eqref{eq:counting-random-variables-no-skip}, 
and $\widehat{P}_{\lambda}$ at \eqref{eq:empirical-transition-matrix-no-skip}.
\end{theorem}

The proof 
(see Section~\ref{section:proof-empirical-estimation-kernel-infinity})
relies on \citet{cohen2020learning} that uses empirical Rademacher averages to derive dimensionless confidence intervals for the problem of learning a distribution with respect to total variation.

\paragraph{Comparison with previous work.}

We note that Theorem~\ref{theorem:fully-empirical-transition-kernel-learning} is superior to \citet[Lemma~9]{pmlr-v117-wolfer20a}. It is instance specific
and $(i)$ only explicitly depends on the number of states by a logarithmic
factor, $(ii)$ can leverage sparsity of the transition matrix (or even more fine grained structure, see e.g. the subsequent examples in this section).
With this approach, we will later show (Theorem~\ref{theorem:learning-s-step-transition-matrix}) how to recover, from a concentration inequality for Markov chains, an upper bound of
\begin{equation}
\label{eq:simplification-learning-kernel-from-empirical}
\frac{c}{\pimin} \max \set{  \frac{1}{\eps^2} \max \set{ \max_{x \in \calX}\nrm{e_x P}_{1/2}, \ln \frac{1}{\delta \eps \pimin} } , \frac{1}{\pssg } \ln \frac{\abs{\calX} \nrm{\mu/ \pi}_{\pi}}{\delta} }    
\end{equation}
for the minimax sample complexity.
As a result, we see that starting from the fully empirical bound of Theorem~\ref{theorem:fully-empirical-transition-kernel-learning} slightly degrades the bound of \citet[Theorem~3.1]{wolfer2021} by a logarithmic factor in the worst case scenario where $\max_{x \in \calX}\nrm{e_x P}_{1/2}$ is of the order of $\abs{\cal{X}}$ and $\log 1 / \pimin > \abs{\calX}$.
However, the new approach becomes vastly superior when $\max_{x \in \calX}\nrm{e_x P}_{1/2}$ can be bounded more tightly (see the examples in this section).
For further comparison, Theorem~\ref{theorem:fully-empirical-transition-kernel-learning} is also stronger than the empirical bound in \citet[Lemma~21]{pmlr-v99-wolfer19a}, which relied on a matrix version of Freedman's inequality and is encumbered with an explicit dependence in $\sqrt{\abs{\calX}}$.

\begin{example}
\label{example:bounded-degree-markov-chains-improves-learning-kernel}
Consider the set of stationary, symmetric and ergodic Markov chains over a state space $\calX$ and with associated connection graph of degree bounded by $\Delta \in \N$, with $\Delta \leq \log \abs{\calX}$.
The stationary distribution of a symmetric chain is uniform, thus $\pimin = \abs{\calX}^{-1}$. Moreover, for any $x \in \calX$, the bounded degree condition is synonymous with the conditional distribution $e_x P$ being sparse. 
With $\calD_x \eqdef \set{ x' \in \calX \colon  P(x,x') \neq 0}$, we have $\abs{\calD_x} \leq \Delta$ and
\begin{equation*}
    \nrm{e_x P}_{1/2}^{1/2} = \sum_{x' \in \calX} \sqrt{P(x,x')} = \sum_{\substack{x' \in \calD_x}} \sqrt{P(x,x')}.
\end{equation*}
It is not hard to see (e.g. with the method of Lagrange multipliers) that the distribution maximizing the above sum is uniformly concentrated over $\calD_x$, and thus
\begin{equation*}
    \nrm{e_x P}_{1/2}^{1/2} \leq \sum_{\substack{x' \in \calD_x}} \frac{1}{\sqrt{\Delta}} \leq \sqrt{\Delta}.
\end{equation*}
As a result,
\eqref{eq:simplification-learning-kernel-from-empirical} can be controlled more simply by
\begin{equation}
c \abs{\calX} \max \set{ \frac{1}{ \eps^2} , \frac{1}{\pssg}  } \ln \frac{\abs{\calX}}{\delta \eps},
\end{equation}
leading to a savings of roughly $\abs{\calX}$.
\end{example}

\begin{example}
\label{example:bounded-renyi-entropy-chains}
Even if the connection graph associated to the Markov chain is not sparse, as long as
the conditional distribution defined for each state is low--entropic, the sample complexity still enjoys a speed--up.
For $\alpha > 0$, and $\mu \in \calP(\calX)$, recall the definition of the R\'{e}nyi entropy,
\begin{equation*}
\begin{split}
H_\alpha(\mu) \eqdef \frac{1}{1 - \alpha} \ln \left( \sum_{x \in \calX} \mu(x)^\alpha \right),
\end{split}
\end{equation*}
which is directly related to the $\ell_\alpha$ quasi--norm of $\mu$. 
Monotonicity properties of $H_\alpha$ with respect to $\alpha$ then implies 
if there exists $\alpha \in (0, 1/2]$ and $\overline{H}_\alpha > 0$ 
such that for any state $x \in \calX$, $H_\alpha(e_x P) < \overline{H}_\alpha$, 
we can further bound
\begin{equation*}
    \max_{x \in \calX}\nrm{e_x P}_{1/2} \leq \exp \left(\frac{1 - \alpha}{\alpha} \overline{H}_\alpha\right).
\end{equation*}
\end{example}

As a consequence of Theorem~\ref{theorem:fully-empirical-transition-kernel-learning} and continuity of $\kappa_\star$ with respect to $\ell_\infty$ \citep[Fact~5.1]{pmlr-v117-wolfer20a}, we readily obtain a confidence interval for $\widehat{\kappa}_\lambda \eqdef \cc (\widehat{P}_{\lambda} )$, the plug--in estimator of the Dobrushin coefficient $\kappa$ of $P$.

\begin{corollary}[Empirical estimation of Dobrushin contraction coefficient $\kappa$]
\label{corollary:estimate-kappa}
Let $P$ be an ergodic transition matrix over a finite space $\calX$, 
$\mu \in \calP(\calX)$, $\delta \in (0,1)$, 
$\lambda \in (0,\infty)^{\calX}$ a smoothing vector, 
$m \in \N$, and $X^m = X_1, X_2, \dots, X_m \sim (\mu, P)$.
With probability at least $1 - \delta$, it holds that
\begin{equation*}
\begin{split}
\abs{\widehat{\kappa}_{\lambda}(X^m) - \kappa} &\leq \widehat{W}_{\lambda, \delta}(X^m),\\
\end{split}
\end{equation*}
with $\widehat{W}_{\lambda, \delta}(X^m)$ defined as in Theorem~\ref{theorem:fully-empirical-transition-kernel-learning}.
\end{corollary}

We conclude this section with an asymptotic analysis of the obtained confidence interval.

\begin{lemma}
\label{lemma:asymptotic-behavior-matrix}
Let $P$ be an ergodic transition matrix over a finite space 
$\calX$ with stationary distribution $\pi$,
$\delta \in (0,1)$, 
$\lambda \in (0 ,\infty)^\calX$ a smoothing vector.
Writing
\begin{equation}
\label{equation:asymptotic-quantity-matrix}
\begin{split}
W^\infty_\delta(P) &\eqdef \frac{2}{\sqrt{m}} \max_{x \in \calX } \set{ \frac{1}{\sqrt{\pi(x)}} \left( \nrm{e_x P}_{1/2}^{1/2} + 3 \sqrt{\frac{1}{2} \ln \frac{2 m \abs{\calX}}{ \delta}} \right)}, \\
\end{split}
\end{equation}
it holds that
\begin{alignat*}{2}
\frac{\widehat{W}_{\lambda, \delta}(X^m)}{W^\infty_\delta(P)} &\as 1, \qquad \sqrt{\frac{m}{\ln{m}}} \widehat{W}_{\lambda, \delta}(X^m) &&\as \frac{3}{\sqrt{2 \pimin}}.
\end{alignat*}
\end{lemma}

\begin{remark}
    Observe that the smoothing vector $\lambda$ does not appear in the definition of $W_\delta^\infty$. In fact, when $\nrm{\lambda}_\infty \leq 1$, the impact of $\lambda$ on the estimation rate is minimal. The primary role of the smoothing vector $\lambda$ is to ensure that all our estimators are well--defined.
\end{remark}

\section{Empirical estimation of \texorpdfstring{$\tgencc$ and $\tmix$}{the generalized contraction coefficient and the mixing time}}
\label{section:empirical-estimation-tmix}
In this section, we design a fully empirical, high--confidence interval for $\tmix$, 
after observing one single long run $X^m = X_1, X_2, \dots, X_m$ from the unknown chain, i.e. for a confidence parameter $\delta \in (0, 1)$, we construct 
an interval $\widehat{T}_{\delta}(X^m) \subset [1, \infty)$, that verifies
\begin{equation*}
\begin{split}
\PR{\tmix \in \widehat{T}_{\delta}(X^m)} &\geq 1 - \delta,
\end{split}
\end{equation*}
with non--trivial $\widehat{T}_{\delta}$.
Our strategy is to derive such a confidence interval for $\tgencc$ instead, defined at \eqref{eq:generalized-contraction-coefficient}, 
and rely on the relationship between $\tmix$ and $\tgencc$ formalized
at Corollary~\ref{corollary:control-tmix-with-kgen-alt2020} to deduce the interval for $\tmix$.
In Section~\ref{section:estimator-definition}, we introduce the estimator for $\tgencc$.
In Section~\ref{section:empirical-confidence-intervals}, 
Lemma~\ref{lemma:explicit-empirical-intervals} yields a confidence 
interval for the estimation of $\tgencc$, Theorem~\ref{theorem:empirical-estimation-tmix} establishes the confidence interval for $\tmix$ together with an asymptotic analysis.
All proofs are deferred to 
Section~\ref{section:proofs} 
for the sake of readability.

\subsection{Estimator definition}
\label{section:estimator-definition}

Recall that the definition of $\tgencc$ in \eqref{eq:generalized-contraction-coefficient}
dauntingly involves taking a maximum over the natural integers.
Our estimator will be the truncated plug--in version of 
$\tgencc$ adapted (smoothed) from \citet{pmlr-v117-wolfer20a}, 
where we only explore a finite prefix $[S]$ of the integers, 
a similar idea to that employed in \cite{pmlr-v99-wolfer19a} for 
estimating the pseudo--spectral gap \eqref{eq:pseudo-spectral-gap}. 
We begin by extending counting random variables in \eqref{eq:counting-random-variables-no-skip} and \eqref{eq:empirical-transition-matrix-no-skip}, to accommodate for a skipping rate $s \in [m - 1]$,
\begin{equation}
\label{eq:counting-random-variables-skip}
\begin{split}
N_{x}^{(s)} &\eqdef \sum_{t=1}^{\floor{(m-1)/s}} \pred{X_{1 + s(t-1)} = x}, \qquad \Nmin^{(s)} \eqdef \min_{x \in \calX} N_x^{(s)} , \\
N_{xx'}^{(s)} &\eqdef \sum_{t=1}^{\floor{(m-1)/s}} \pred{X_{1 + s(t-1)} = x, X_{1 + st} = x'},\\
\widehat{P}^{(s)}_{\lambda} &\eqdef \sum_{x,x' \in \calX} \frac{N_{xx'}^{(s)} + \lambda(x) }{N_x^{(s)} + \lambda(x) \abs{\calX}} \; e_x \trn e_{x'}.
\end{split}
\end{equation}
For a fixed $S \in \N$, we then let the estimator for $\tgencc$ parametrized by $S$ be defined by
\begin{equation}
\label{eq:definition-fixed-truncation-estimator}
\esttgencc \colon \calX^m \to (0, 1), \; X^m \mapsto 1 - \max_{s \in [S]} \set{\frac{1 - \estcc_s(X^m)}{s}},
\end{equation}
where for $s \in \N$, we define \footnote{When there is no ambiguity, we henceforth may omit the subscript $\lambda$ in order to simplify notation.}
\begin{equation*}
    \estcc_s(X^m) \eqdef \cc\left(\widehat{P}^{(s)}_{\lambda}\right).
\end{equation*}
In words, the estimator will compute the empirical contraction coefficients for the first $S$ skipping rates,
and return a value that involves maximizing a ratio over this finite set.
For any fixed $S$, we know from Theorem~\ref{theorem:arbitrary-large-s} that we can anticipate an approximation error.
Let us highlight the following trade--off. 
The larger $S$, the smaller we expect our approximation error to be. In fact, writing $$\tgenccS \eqdef 1- \max_{s \in [S]} \set{ (1 - \cc_s)/s },$$ 
it holds that $\abs{\tgencc - \tgenccS} \leq 1/S$.
However, to avoid unnecessary computation,
we would prefer to take a small value for $S$,
which still guarantees convergence.
We devise a fully--adaptive method that lets also $S$ vary with the sample,
and deduce the following estimator for $\tgencc$
\begin{equation}
\label{eq:definition-general-estimator-adaptive}
\esttgenccadapt \eqdef \esttgenccexplest,
\end{equation}
where $\widehat{S} \colon \calX^m \to \N$ is defined by
\begin{equation}
\begin{split}
\label{eq:adaptive-prefix}
\widehat{S}(X^m) &\eqdef \ceil*{\frac{1}{ \widehat{K}^\phi_{\lambda, \delta}(X^m)}}, \qquad 
\widehat{K}^\phi_{\lambda, \delta}(X^m) \eqdef \max_{s \in [\phi(m)]} \set{ \frac{1}{s} \widehat{W}^{(s)}_{\lambda, \delta}(X^m) }
\end{split}
\end{equation}
where for any $s \in [S]$,
\begin{equation}
\label{eq:empirical-interval-skipped-chain}
\begin{split}
 \widehat{W}_{\lambda, \delta}^{(s)}(X^m) \eqdef \widehat{W}_{\lambda, \delta}\left(X^m_{[s]}\right)
\end{split}
\end{equation}
is the empirical interval defined in Theorem~\ref{theorem:fully-empirical-transition-kernel-learning} 
computed for the $s$--skipped Markov chain,
and $\phi \colon \N \to \N$ is a non--decreasing function such that $\lim_{m \to \infty} \phi(m) = +\infty$. For instance, it is appropriate to choose a slowly increasing function, such as a logarithm.

\subsection{Empirical confidence intervals}
\label{section:empirical-confidence-intervals}

We now derive a confidence interval for the problem of estimating $\tgencc$,
the generalized contraction coefficient of the chain.

\begin{lemma}
\label{lemma:explicit-empirical-intervals}
Let $P$ be an ergodic transition matrix over a finite space $\calX$, 
with generalized contraction coefficient $\tgencc$ \eqref{eq:generalized-contraction-coefficient},
$\mu \in \calP(\calX)$, $\delta \in (0,1)$, 
$\lambda \in (0, \infty)^{\calX}$ a smoothing vector, 
$m \in \N$, and $X^m = X_1, X_2, \dots, X_m \sim (\mu, P)$.
\begin{enumerate}
	\item[$(a)$] Let $S \in \N$. With probability at least $1 - \delta$, it holds that
\begin{equation*}
\begin{split}
\abs{\esttgencc - \tgencc} &\leq \widehat{K}_{\lambda, \delta, S}(X^m), \\
\end{split}
\end{equation*}
with
\begin{equation*}
\begin{split}
\widehat{K}_{\lambda, \delta, S}(X^m) &\eqdef \frac{1}{S} + \max_{s \in [S]} \set{\frac{1}{s} \widehat{W}_{\lambda, \delta/S}^{(s)}(X^m)},
\end{split}
\end{equation*}
and where $\esttgencc$ is the fixed truncation estimator defined at \eqref{eq:definition-fixed-truncation-estimator}, 
and $\widehat{W}_{\lambda, \delta}^{(s)}(X^m)$ is defined in \eqref{eq:empirical-interval-skipped-chain}.
\item[$(b)$] 
With probability at least $1 - \delta$, it holds that
\begin{equation*}
\begin{split}
\abs{\esttgenccadapt - \tgencc} &\leq \widehat{K}_{\lambda, \delta}(X^m), \\
\end{split}
\end{equation*}
with
\begin{equation*}
\begin{split}
\widehat{K}_{\lambda, \delta}(X^m) &\eqdef \frac{1}{\widehat{S}} + \max_{s \in [\widehat{S}]} \set{\frac{1}{s} \widehat{W}_{\lambda, \delta/( \ceil{\sqrt{m}}\widehat{S})}^{(s)}(X^m)},
\end{split}
\end{equation*}
and where the adaptive truncation $\widehat{S}$ and estimator $\esttgenccadapt$ are defined respectively at \eqref{eq:adaptive-prefix} and \eqref{eq:definition-general-estimator-adaptive}.
\end{enumerate}
\end{lemma}

\paragraph{Comparison with previous work.}

\begin{enumerate}[label=$(\roman*)$]
    \item Lemma~\ref{lemma:explicit-empirical-intervals}--$(a)$ enjoys the same instance dependent properties as Lemma~\ref{theorem:fully-empirical-transition-kernel-learning}, thus compares favorably with \citet[Theorem~2]{pmlr-v117-wolfer20a}. 
    \item The confidence interval computed in Lemma~\ref{lemma:explicit-empirical-intervals}--$(a)$ is narrower than the one designed around the pseudo--spectral gap in \cite{wolfer2023improved}. Even taking the unfavorable case where $\max_{x \in \calX}\nrm{e_x \widehat{P}_\lambda}_{1/2}$ is of the order of $\abs{\calX}$, Lemma--$(a)$~\ref{lemma:explicit-empirical-intervals} yields (ignoring logarithmic terms),
    \begin{equation}
    \label{eq:asymp-decay-k-gen-truncated}
\begin{split}
\abs{\esttgencc - \tgencc} - \frac{1}{S} \asymp \sqrt{\frac{\abs{\calX}}{\pimin m}},
\end{split}
\end{equation}
whereas the known intervals around the pseudo--spectral gap $\pssg$ of the estimator $\estpssgS$ 
defined by \citeauthor{wolfer2023improved} could generally decay as slowly as
\begin{equation}
\begin{split}
\label{eq:asymp-rates-pssg-truncated}
\abs{\estpssgS - \pssg} - \frac{1}{S} \asymp \sqrt{\frac{\log^3 m}{m}} \frac{\log 1/\pi_\star}{\pssg \pi_\star^{3/2}}.
\end{split}
\end{equation}
    \item Notably, the asymptotic rate in \eqref{eq:asymp-decay-k-gen-truncated} does not depend on the mixing properties of the unknown chain. This stands in contrast with the analysis performed in both \citet{hsu2019} and \citet{wolfer2023improved}, that introduces the inverse of the (pseudo--)spectral gap in the rate \eqref{eq:asymp-rates-pssg-truncated}, by relying on a perturbation analysis of the stationary distribution.
    \item The confidence interval is also easier to compute from the data than in \citet{hsu2019}. In particular, it does not involve the computation of a group (Drazin) inverse \citep{meyer1975role}.
    \item Lemma~\ref{lemma:explicit-empirical-intervals}--$(b)$ proposes a fully adaptive method for computing $\tgencc$, instead
of a hard--coded number of skipping rates to consider together with some heuristic value \citep[(14)]{pmlr-v117-wolfer20a}. The prefix of integers to explore defined at \eqref{eq:adaptive-prefix} is set explicitly.
\end{enumerate}

As a direct consequence of Lemma~\ref{lemma:explicit-empirical-intervals}, 
and since from Corollary~\ref{corollary:control-tmix-with-kgen-alt2020}, 
$\tmix$ and $1/(1 - \tgencc)$ are within universal multiplicative constants,
we deduce an empirical confidence interval for the mixing time itself.

\begin{theorem}
\label{theorem:empirical-estimation-tmix}
Let $P$ be an ergodic transition matrix over a finite space 
$\calX$ with mixing time $\tmix$, 
$\mu \in \calP(\calX)$, $\delta \in (0,1)$, 
$\lambda \in (0, \infty)^\calX$ a smoothing vector, 
$m \in \N$, and $X^m = X_1, X_2, \dots, X_m \sim (\mu, P)$.
With probability at least $1 - \delta$, it holds that
$\tmix \in \widehat{T}_{\lambda, \delta}(X^m)$, with
\begin{equation*}
\begin{split}
\widehat{T}_{\lambda, \delta}(X^m) \eqdef \left(  \frac{\underline{c}}{1 - \floor{\esttgenccadapt(X^m) - \widehat{K}_{\lambda, \delta/2}(X^m)}_+}, \frac{\overline{c}}{\floor{1 - \esttgenccadapt(X^m) - \widehat{K}_{\lambda, \delta/2}(X^m)}_+}\right),
\end{split}
\end{equation*}
and where $\underline{c}, \overline{c}$ are universal constants defined in Corollary~\ref{corollary:control-tmix-with-kgen-alt2020},
$\esttgenccadapt(X^m)$ is the estimator for $\tgencc$ defined at \eqref{eq:definition-general-estimator-adaptive}, and
$\widehat{K}_{\lambda, \delta/2}(X)$ is defined in Lemma~\ref{lemma:explicit-empirical-intervals}.
\end{theorem}

We conclude this section on empirical estimation with an asymptotic width 
analysis of the empirical confidence interval.

\begin{theorem}
\label{theorem:asymptotic-behavior-mixing-parameters}
Let $P$ be an ergodic transition matrix over a finite space 
$\calX$ with stationary distribution $\pi$,
$\delta \in (0,1)$, 
$\lambda \in (0 ,\infty)^\calX$ a smoothing vector and $S \in \N$.

\begin{equation*}
    \sqrt{\frac{m}{\log m} } \left(\widehat{K}_{\lambda, \delta, S}(X) - \frac{1}{S} \right) \as \frac{3}{\sqrt{2 \pimin}}.
\end{equation*}

\end{theorem}

\section{Point estimators}
\label{section:minimax-estimation}

In this section, we depart from fully empirical estimation and move towards characterizing the sample complexity of estimating the mixing time to multiplicative accuracy, as opposed to one of its spectral proxies ($\gamma_\star, \pssg$). As in previous works, our upper bounds will necessarily involve oracle quantities, in particular the mixing time itself.

\subsection{Estimation of \texorpdfstring{$P^s$}{powers of transition kernels}}
\label{section:minimax-estimation-kernel}
In this section, we estimate the power of a transition kernel with a skipped chain, by considering the counting random variables introduced in the previous section \eqref{eq:counting-random-variables-skip}.

\begin{theorem}
\label{theorem:learning-s-step-transition-matrix}
Let $\eps, \delta \in (0,1), m \in \N$, $s \in [m-1]$, and a finite state space $\calX$. 
There exists an estimator 
$\widehat{P}^{(s)}_{\lambda} \colon \calX^m \to \calW (\calX)$,
such that for any ergodic transition matrix $P$ over 
$\calX$, $\mu \in \calP(\calX)$,
$X^m = X_1, X_2, \dots, X_m \sim (\mu, P)$,
and 
$$m \geq c \max \set{ \frac{s \Gamma(P^s)}{\eps^2},  \frac{s \nrm{\lambda/\pi}_{\infty} \abs{\calX}}{\eps}, \frac{s}{\eps^2 \pimin} \ln \frac{1}{\delta \eps \pimin}, \frac{1}{(1 - 
\kappa_\star) \pimin} \ln \frac{\abs{\calX} \nrm{\mu/ \pi}_{\pi}}{\delta} }$$
it holds\footnote{We can obtain a similar bound with $\pssg$ in lieu of $1 - \kappa_\star$.} that
$$\nrm{\widehat{P}^{(s)}_{\lambda} - P^s}_{\infty} \leq \eps,$$
with probability at least $1 - \delta$, where $\pimin$ and $\kappa_\star$ pertain to $P$, $\nrm{\mu/\pi}_{\pi}$ is defined in \eqref{eq:pi-norm}, and
\begin{equation}
\label{eq:definition-useful-chain-complexities}
\begin{split}
\nrm{\lambda/\pi}_{\infty} &\eqdef \max_{x \in \calX} \set{\frac{\lambda(x)}{\pi(x)}}, \qquad
\Gamma(P) \eqdef \max_{x \in \calX} \set{ \frac{\nrm{e_x P}_{1/2}}{\pi(x)}} \leq \frac{\abs{\calX}}{\pimin}.
\end{split}
\end{equation}
\end{theorem}
In addition to generalizing \citet[Theorem~3.1]{wolfer2021} to the problem of learning an arbitrary power $s \in [m-1]$ of the transition operator, Theorem~\ref{theorem:learning-s-step-transition-matrix} refines the result, even for $s = 1$, allowing for leveraging additional information about the kernel, as $\Gamma(P) \leq \abs{\calX} / \pimin$ (see the example below).

The quantity $\Gamma(P)$ introduces an interesting trade--off, between $(a)$ the half--norm of the conditional distribution of a particular state, which measures exactly \citep{cohen2020learning} its hardness for the learning problem, and $(b)$ its stationary probability, that governs its long--run appearance in the trajectory.

\begin{example}
Consider $U$, the uniform transition matrix encoding the transition probabilities of the max--entropic chain over the fully connected graph $(\calX, \calX^2)$. Immediately, $\Gamma(U) = \abs{\calX}^2$.
On the other hand, let $P$ with an associated connection graph that is $2$--regular. We have $\abs{\calX} \leq \frac{1}{\pimin} \leq \Gamma(P) \leq \frac{2}{\pimin}$. Depending on the value of $\pimin$, $\Gamma(P)$ can be much smaller or much larger than $\Gamma(U)$.
\end{example}

\subsection{Estimation of \texorpdfstring{$\tmix$}{the mixing time} with multiplicative accuracy}
\label{section:minimax-estimation-mixing-properties}
For chosen precision $\eps$ and confidence $\delta$ parameters, we first construct a point estimator for $
\kappa_\star$, down to additive error, where the algorithm only needs knowledge of $\abs{\calX}, \eps$ and $\delta$ in order to run. To reach arbitrary precision, we explore the first 
$S = \ceil{2/\eps}$ possible skipping rates, i.e. consider the estimator $\widehat{\kappa}_\star^+ \eqdef  \esttgenccexpl$ \eqref{eq:definition-fixed-truncation-estimator}. 

\begin{theorem}[Point estimator for $\tgencc$ in additive error]
\label{theorem:learn-tgencc-ub-absolute}
Let $\eps, \delta \in (0,1)$ and $m \in \N$. 
There exists an estimation procedure 
$\esttgenccpea \colon \calX^m \to (0,1)$ such that 
for any unknown ergodic Markov chain with transition matrix $P$, minimum stationary probability
$\pimin$ and generalized contraction coefficient $\tgencc$,
$X^m = X_1, X_2, \dots, X_m \sim (\mu, P)$, if 
\begin{equation*}
m \geq  c \max \set{ \frac{1}{\eps^2}\max_{s \leq \ceil{2/\eps}} \set{ \frac{\Gamma(P^s)}{s} }, \frac{1}{\eps^2 \pimin} \ln \frac{1}{\delta \eps \pimin}, \frac{1}{(1 - \tgencc)\pimin}\ln \frac{1}{\pimin \delta} },
\end{equation*}
then $\abs{\esttgenccpea - \tgencc} < \eps$ holds with probability at least $1-\delta$, where $c > 0$ is a universal constant, and $\Gamma$ is defined at \eqref{eq:definition-useful-chain-complexities}.
\end{theorem}

\begin{example}
\label{example:bounded-degree-markov-chains-improves-learning-kappa-absolute}
Recall the set of symmetric and ergodic Markov chains with $\Delta$--bounded degree.
When the state space is large and $\Delta \ll \abs{\calX}$, we can still obtain 
$$\max_{s \in \ceil{2/\eps}} \set{\Gamma(P^s)/s} \ll \abs{\calX}^2.$$
We note that
for any $s \in \N$, and any $x \in \calX$, 
$$\max_{x \in \calX} \nrm{e_x P^s}_{1/2}^{1/2} \leq \left( \max_{x \in \calX}\nrm{e_x P}_{1/2}^{1/2}\right)^s.$$
To see this, let $x_0 \in \calX$. For $s = 1$, the claim is trivial. 
Let now $s > 1 $. By definition,
\begin{equation*}
\begin{split}
\nrm{e_{x_0} P^s}_{1/2}^{1/2} &= \sum_{x' \in \calX} \sqrt{P^s(x_0,x')} = \sum_{x' \in \calX} \sqrt{\sum_{x'' \in \calX} P^{s-1}(x_0, x'') P(x'',x')} \\
&\stackrel{(i)}{\leq} \sum_{x' \in \calX} \sum_{x'' \in \calX} \sqrt{P^{s-1}(x_0, x'') P(x'',x')} \\
&= \sum_{x'' \in \calX} \sqrt{P^{s-1}(x_0, x'') } \sum_{x' \in \calX}\sqrt{P(x'',x')} \\
&\leq \nrm{e_{x_0} P^{s-1}}_{1/2}^{1/2} \max_{x \in \calX}\nrm{e_x P}_{1/2}^{1/2} \stackrel{(ii)}{\leq} \left( \max_{x \in \calX}\nrm{e_x P}_{1/2}^{1/2}\right)^s, \\
\end{split}
\end{equation*}
where $(i)$ follows from $a,b > 0 \implies \sqrt{a + b} \leq \sqrt{a} + \sqrt{b}$, 
and $(ii)$ is by induction and taking a maximum over the state space.
From there,
\begin{equation*}
    \max_{s \in \ceil{2/\eps}} \set{\frac{\Gamma(P^s)}{s}} \leq  \frac{\abs{\calX} \Delta^{\ceil{2/\eps}}}{\ceil{2/\eps}}.
\end{equation*}
\end{example}

We proceed to construct an algorithm that outputs an estimate 
of $1 - \tgencc$ with constant multiplicative error.
The solution is an amplification argument that is
fleshed out in
Section~\ref{section:proof-point-estimator-relative-error}. 
We consider the estimator 
\begin{equation}
\label{eq:amplified-estimator-def}
\widehat{\kappa}^\times_\star \eqdef 1 - \frac{1 - \widehat{\kappa}^{(S)}}{\sqrt{c_\star} S},
\end{equation}
where the natural constant $c_\star$ appears in Proposition~\ref{proposition:kgen-of-power},
$S$ is a random skipping rate defined by
\begin{equation*}
    \log_2 S \eqdef \argmin_{p \in \set{0,1, 2, \dots}} \set{ 1 - \widehat{\kappa}_\star^{(2^p)} > \frac{\sqrt{2}}{4}},
\end{equation*} and where for $s \in \N$, $\widehat{\kappa}_\star^{(s)}$ is the estimator for $\kappa_\star(P^s)$ defined in Theorem~\ref{theorem:learn-tgencc-ub-absolute} that achieves additive precision $\eps$ with the skipped chain $X_{[s]}^m$ as input.

\begin{theorem}[Point estimator for $1 - \tgencc$ in multiplicative error]
\label{theorem:learn-tgencc-ub-relative}
Let $\delta \in (0,1)$, and $m \in \N$. 
The estimation procedure 
$\widehat{\kappa}_\star^{\times} \colon \calX^m \to (0,1)$ defined in \eqref{eq:amplified-estimator-def} is such that
for any unknown ergodic Markov chain with transition matrix $P$, minimum stationary probability
$\pimin$ and coefficient $\tgencc$,
$X^m = X_1, X_2, \dots, X_m \sim (\mu, P)$, if 
\begin{equation*}
m \geq  c \max \set{ \Xi(P),  \frac{1}{(1 - \kappa_\star)\pimin}\ln \frac{\log (1 - \kappa_\star)^{-1}}{\delta \pimin } },
\end{equation*}
then 
\begin{equation*}
    \max \set{ \frac{1 - \widehat{\kappa}^\times_\star}{1 - \kappa_\star}, \frac{1 - \kappa_\star}{1 - \widehat{\kappa}^\times_\star}} \leq 2\sqrt{2\log 4e} = 4.36925\dots
\end{equation*}
holds with probability at least $1-\delta$, where
\begin{equation}
\label{eq:xi-measure-complexity}
    \Xi(P) \eqdef \max_{p \in \set{0,1, \dots, \overline{p}}} \set{ 2^p \max_{s \leq 16} \set{ \frac{\Gamma\left(P^{s 2^p}\right)}{s}}} \leq  \frac{\abs{\calX}}{(1 - \kappa_\star)\pimin},
\end{equation}
$c > 0$ is a universal constant, $\overline{p} \eqdef \ceil{\log_2 \tmix}$, and $\Gamma$ is defined at \eqref{eq:definition-useful-chain-complexities}.
\end{theorem}

\begin{remark}
The doubling trick is inspired by the one of \citet{levin2016estimating, hsu2019} in the context of estimating the spectral gap. However, the latter is not applicable verbatim to our problem. The reason is that whereas it is possible to relate exactly the absolute spectral gap of a Markov chain and its $s$--skipped version by the following equality
\begin{equation*}
\gamma_\star(P^s) = 1 - (1 - \gamma_\star(P))^s,
\end{equation*}
the correspondence between $\kappa_\star(P)$ and $\kappa_\star(P^s)$ is not as straightforward (see Proposition~\ref{proposition:kgen-of-power}). As a consequence, we are unable to produce a bound to arbitrary multiplicative precision, and must settle for a universal constant. 
\end{remark}

\begin{remark}
For any $a, b > 0$, it holds that
$\abs{a/b - 1} \leq \max \set{a/b, b/a} - 1$. Consequently, for a trajectory length upper bounded as in Theorem~\ref{theorem:learn-tgencc-ub-relative}, we can rewrite that
$$\abs{ \frac{1 - \widehat{\kappa}^\times_\star}{1 - \kappa_\star} - 1 } \leq 2\sqrt{2\log 4e} - 1$$
holds with probability at least $1-\delta$.
\end{remark}

Finally, we convert our amplified estimator of the contraction coefficient into an estimator for the mixing time,
\begin{equation}
\label{eq:esttmix-def}
    \esttmix(X^m) \eqdef \frac{\sqrt{\overline{c} \underline{c}}}{1 - \widehat{\kappa}_\star^{\times}(X^m)},
\end{equation}
where the constants $\underline{c}$ and $\overline{c}$ are defined in Corollary~\ref{corollary:control-tmix-with-kgen-alt2020}.

\begin{corollary}[to Theorem~\ref{theorem:learn-tgencc-ub-relative}, in terms of mixing time]
\label{corollary:estimator-tmix-relative}
Let $\delta \in (0,1)$ and $m \in \N$. 
The estimator 
$\esttmix \colon \calX^m \to \N$ defined in \eqref{eq:esttmix-def} is such that 
for any unknown ergodic Markov chain with transition matrix $P$, minimum stationary probability
$\pimin$ and mixing time $\tmix$,
$X^m = X_1, X_2, \dots, X_m \sim (\mu, P)$, if 
\begin{equation*}
m \geq  c \max \set{ \Xi(P),  \frac{\tmix}{\pimin}\ln \frac{\log \tmix }{\delta \pimin } },
\end{equation*}
then 
\begin{equation*}
\max \set{ \frac{\esttmix}{\tmix}, \frac{\tmix}{\esttmix}} \leq \sqrt{\frac{\overline{c}}{\underline{c}}} 4\sqrt{\log 4e}
\end{equation*}
holds with probability at least $1-\delta$, where $c > 0$ is a universal constant, and $\Xi(P)$ is as in
Theorem~\ref{theorem:learn-tgencc-ub-relative}.
\end{corollary}

\paragraph{Comparison with previous work.}

In order to estimate the relaxation time of a Markov chain to multiplicative accuracy $\eps$, it suffices ---at least in the reversible setting \citep{levin2016estimating, hsu2019}--- to observe a trajectory of length about
\begin{equation}
\label{eq:ub-asg-rel}
    m \asymp \frac{\trel}{\pimin \eps^2}.
\end{equation}
For comparison, considering the worst--case where $\Xi(P)$ and $\tmix \abs{\calX}/\pimin$ are of the same order,
Theorem~\ref{theorem:learn-tgencc-ub-relative} yields an upper bound of
\begin{equation}
\label{eq:ub-kgen-rel}
m \asymp  \frac{\tmix \abs{\calX}}{\pimin},
\end{equation}
for the problem of estimating $\tmix$ to constant multiplicative precision.
    We further note that \eqref{eq:ub-asg-rel} holds for arbitrary multiplicative accuracy $\eps$, while our guarantee is for constant precision.
    \citet[Theorem~6]{pmlr-v117-wolfer20a}, through a different approach, obtained a sample complexity upper bound of
    \begin{equation*}
        m \asymp \frac{\tmix^2 \abs{\calX}}{\pimin \eps^2},
    \end{equation*}
    which can achieve the arbitrary multiplicative precision $\eps$, but carries an extra multiplicative factor in $\tmix$.
    Nevertheless, the logarithmic gap between $\trel$ and $\tmix$ precludes estimates of $\trel$ to yield corresponding estimates to arbitrary accuracy for $\tmix$, thus our method eventually estimates $\tmix$ more closely.

\section{Analysis of a family of three state chains}
\label{section:family-special}
In this section, we investigate additional properties of the generalized contraction coefficient defined in \eqref{eq:generalized-contraction-coefficient}. Recall that computing $\tgencc$ requires a search over the integers that correspond to the different skipping rates one may consider.
If the chain is over two states, then $\tgencc$ is completely characterized
(see Lemma~\ref{lemma:properties-two-state-mc}), 
and the smallest integer that achieves $\tgencc$ is $s = 1$.
For state spaces of larger sizes, one may wonder whether there exists an upper bound on the integer $s$ that depends solely on the state space size $\abs{\calX}$.
For comparison, from Wielandt's theorem (see e.g. \citet{schneider2002wielandt} for a proof), 
for a stochastic matrix $P$, a sharp bound for the exponent of primitivity \footnote{The smallest $p \in
\N$ such that $A^p > 0$ entry-wise.} exists, and is given by $(\abs{\calX}^2 - 1)^2 + 1$.
However here, we show that even if the chain has only three states, the exact solution of this optimization problem in the definition \eqref{eq:generalized-contraction-coefficient} can require us to look at arbitrary large integers.
To demonstrate this fact, we introduce the family of chains parametrized by $\eta \in (0, 2/3)$,
$$P_\eta = \begin{pmatrix} 1/3 & 0 & 2/3 \\ \eta/2 & 1 - \eta & \eta/2 \\ \eta/2 & \eta/2 & 1 - \eta \end{pmatrix},$$
with stationary distribution given by
$$\pi_\eta = \frac{8}{3(3 \eta + 4)} \left(\frac{9 \eta}{8}, 1/2, 1 \right),$$
and exhibit some of its properties.

\begin{lemma}
\label{lemma:kgen-special-family}
For $\eta \in (0, 1/6)$, and $s \in \N$,
\begin{enumerate}
	\item[$(a)$] $\cc_s(P_\eta) = \left( 1 - \frac{3 \eta}{2} \right)^{s -1} \left( 1 - \left( \frac{1 + 3^{-s}}{2} \right) \frac{3 \eta}{2} \right).$
    \item[$(b)$] $\cc_1(P_\eta) = 1 - \eta$.
    \item[$(c)$] $(1 - \cc_2(P_\eta))/2 > 1 - \cc_1(P_\eta)$.
\end{enumerate}
\end{lemma}

Lemma~\ref{lemma:kgen-special-family} shows that 
the chain induced from $P_\eta$ is contractive,
and that it suffices to take $\eta < 4/15$ for the maximum $\max_{s \in N} \set{ (1 - \cc_s) /s}$ to not be achieved for $s = 1$. This also yields the somewhat simple expression for $\tgencc$,
\begin{equation*}
\begin{split}
\tgencc(P_\eta) &= 1 - \max_{s \in \N} \set{\theta(s, \eta)}, \\ \theta(s, \eta) &\eqdef  \frac{1}{s} \left( 1 - \left( 1 - \frac{3 \eta}{2} \right)^{s -1} \left( 1 - \left( \frac{1 + 3^{-s}}{2} \right) \frac{3 \eta}{2} \right)\right).
\end{split}    
\end{equation*}
This family of chain possesses the surprising property that any $s_\star (\eta)$ that maximizes $\theta(s, \eta)$ increases to infinity as $\eta \to 0$.

\begin{lemma}
\label{lemma:special-category-chain-requires-arbitray-large-s}
For any $s_0 \in \N$, there exists $\eta \in (0, 1/6)$ such that 
$$\argmax_{s \in \N} \set{ \frac{1 - \cc_s(P_\eta)}{s}} \subset [s_0, \infty).$$
\end{lemma}
This construction leads to the general statement that the skipping rate that achieves $\tgencc$ can be arbitrary large, making the search over $\N$ generally unavoidable.

\section{Conclusion, discussion, limitations and further work}
\label{section:conclusion-discussion}
\citet{pmlr-v117-wolfer20a} has closed the logarithmic gap in the mixing time approximation problem by proposing the coefficient $\kappa_\star$, and in this paper, have performed a tighter analysis for estimating $\kappa_\star$ from the data.
It is instructive to compare the minimax upper bound on the rate
for estimating $\trel$ in \eqref{eq:ub-asg-rel}, with the one we obtain for estimating $\tmix$, given roughly by
\begin{equation}
m \asymp  \frac{\tmix}{\pimin} + \Xi(P),
\end{equation}
where the term $\Xi(P)$ depends finely on $P$, but can grow as large as $\tmix \abs{\calX} / \pimin$ in the worst--case. A natural question is the necessity of the additional complexity term.
We leave the formal construction of a corresponding lower bound as an interesting open question.

Furthermore, we have shown that our sample complexity upper bounds for several intermediary problems can benefit from speed--up in favorable cases (e.g. sparsity of the connection graph).
However, for estimating $1- \kappa_\star$ to multiplicative accuracy: $(i)$ it is not clear how the quantity $\Xi(P)$ can be bounded in practice, and $(ii)$ $\Xi(P)$  and $\tmix/\pimin$ are both super--linear in the number of states. The rates are thus prohibitively large for real world MCMC applications that deal with vast state spaces, and our work remains of mostly theoretical interest for this community.

\section{Proofs}
\label{section:proofs}
\subsection{Proof of Theorem~\ref{theorem:control-tau-with-kgen-and-xi-converse}}
Let $\alpha, \beta \in (0, 1)$, $\alpha \neq \beta$, and consider the transition matrix
$P_{\alpha, \beta} = \begin{psmallmatrix} \alpha & 1 - \alpha \\ \beta & 1 - \beta \end{psmallmatrix}$.
From Lemma~\ref{lemma:properties-two-state-mc}, the properties of this chain verify
\begin{equation*}
\begin{split}
\tau(\xi)(1 - \tgencc) &= \ceil*{\frac{\ln{\xi}}{\ln{\abs{\alpha - \beta}}}}(1 - \abs{\alpha - \beta}) \geq \frac{\abs{\alpha - \beta} - 1 }{\ln{\abs{\alpha - \beta}}} \ln{\frac{1}{\xi}} \geq \ln{\frac{1}{\xi}}, \\
\end{split}
\end{equation*}
which entails $(b)$.
Using the same transition matrix, we also show $(c)$ that the lower bound on $\tau(\xi)$ is tight.
Indeed, for $P_{\alpha, \beta}$, whenever $\abs{\alpha - \beta} = \xi$, it holds that $ \tau(\xi)(1 - \tgencc) = 
\ceil*{\frac{\ln{\xi}}{\ln{\abs{\alpha - \beta}}}}(1 - \abs{\alpha - \beta}) = 1 - \xi$ (Lemma~\ref{lemma:properties-two-state-mc}).
Finally, for $n \geq 3$, taking a sequence $P_{\alpha_n, \beta_n}$ such that $\abs{\alpha_n - \beta_n} = \theta_n = \xi + 1/n \in (0, 1)$,
\begin{equation*}
\begin{split}
\lim_{n \to \infty} \ceil*{\frac{\ln{\xi}}{\ln{\theta_n}}}(1 - \theta_n) = 2(1 - \xi),
\end{split}
\end{equation*}
hence $c$ has to verify for any $\xi$, 
$\ln \frac{c e}{\xi} \geq 2(1 - \xi)$, i.e. $c \geq \sup_{\xi \in (0, 1/2)}\frac{\xi}{e} e^{2(1 - \xi)} = 1/2$,
which concludes $(a)$.

\qed

\begin{lemma}
\label{lemma:properties-two-state-mc}
Let $\alpha, \beta \in (0, 1)$, and let
$P_{\alpha, \beta} = \begin{psmallmatrix} \alpha & 1 - \alpha \\ \beta & 1 - \beta \end{psmallmatrix}$. Then,
\begin{equation*}
\begin{split}
	(i) & \qquad \forall s \in \N, \cc_s = \abs{\alpha - \beta}^s, \\
	(ii) & \qquad \tgencc = \abs{\alpha - \beta}, \\
	(iii) & \qquad \forall \xi \in (0, 1/2), \tau(\xi) = \ceil*{\frac{\ln{\xi}}{\ln{\abs{\alpha - \beta}}}}.
\end{split}
\end{equation*}
\end{lemma}
\begin{proof}
A direct computation yields that
\begin{equation*}
\begin{split}
P_{\alpha, \beta}^s &= \frac{1}{1 - \alpha + \beta} \begin{pmatrix} 1 & \frac{\alpha - 1}{\beta} \\ 1 & 1 \end{pmatrix} \begin{pmatrix} 1 & 0 \\ 0 & (\alpha - \beta)^s \end{pmatrix} \begin{pmatrix} \beta & 1 - \alpha \\ -\beta & \beta \end{pmatrix} \\
&= \frac{1}{1 - \alpha + \beta} \begin{pmatrix} \beta - (\alpha - 1)(\alpha - \beta)^s & 1 - \alpha + (\alpha - 1)(\alpha - \beta)^s \\
\beta - (\alpha - \beta)^s \beta & 1 - \alpha + \beta (\alpha - \beta)^s \end{pmatrix},
\end{split}
\end{equation*}
and as a result,
\begin{equation*}
\begin{split}
2 \cc_s &= \frac{2}{1 - \alpha + \beta} \abs{ (\alpha - \beta)^s(\beta - \alpha + 1) } = 2 \abs{\alpha - \beta}^s,
\end{split}
\end{equation*}
which proves claim $(i)$. For $(ii)$, notice that for any $s \in \N$,
\begin{equation*}
\begin{split}
\frac{1 - \abs{\alpha - \beta}^s}{s} &= \frac{1}{s} (1 - \abs{\alpha - \beta}) \sum_{t=0}^{s-1} {\underbrace{\abs{\alpha - \beta}}_{< 1}}^t  \leq 1 - \abs{\alpha - \beta}, \\
\end{split}
\end{equation*}
which is achieved for $s = 1$. Finally, $(iii)$ follows directly from the definition of $\tau(\xi)$ together with $(i)$,
$\tau(\xi) = \min \set{s \in \N \colon \abs{\alpha - \beta}^s < \xi}$, and solving for $s$.
\end{proof}

\subsection{Proof of Proposition~\ref{proposition:kgen-of-power}}
Let $p \in \N$. 
Let $X_1, \dots, X_m \sim (\mu, P)$ a stationary ergodic Markov chain with stationary distribution $\pi$ 
and mixing time $\tmix$. We first verify that the mixing time 
$\tmix^{(p)}$ of the skipped chain,
\begin{equation*}
\begin{split}
X_{1}, X_{1 + s}, X_{1 + 2s}, \dots, X_{1 + \floor{(m-1)/s}s} \sim (\mu, P^s),
\end{split}
\end{equation*}
is such that 
\begin{equation}
\label{equation:tmix-of-power}
\tmix^{(p)} \leq \ceil{\tmix/p}.
\end{equation}
Indeed \citep[Fact~5.2]{pmlr-v117-wolfer20a}, let $t \in \N$ such that $t > \ceil{\tmix/p}$, then
\begin{equation*}
\begin{split}
\tv{\mu (P^p)^t - \pi} &\leq \tv{\mu (P^p)^{\ceil{\tmix/p}} - \pi} \leq \tv{\mu P^{\tmix} - \pi} \leq \frac{1}{4}, \\
\end{split}
\end{equation*}
where the first two inequalities holds as advancing the chain can only move it closer to stationarity 
\citep[Exercise~4.2]{levin2009markov}.
By definition, 
\begin{equation*}
\begin{split}
\tgencc^{(p)} &= 1 - \max_{s \in \N} \set{ \frac{1 - \cc_{sp}}{s} } = 1 - p\max_{s \in \N} \set{ \frac{1 - \cc_{sp}}{sp} } \\
&\geq 1 - p\max_{r \in \N} \set{ \frac{1 - \cc_{r}}{r} } = 1 - p(1 - \tgencc), \\
\end{split}
\end{equation*}
where the inequality follows by taking the maximum over a larger set.
Moreover, from \eqref{equation:tmix-of-power} and Corollary~\ref{corollary:control-tmix-with-kgen-alt2020},
\begin{equation*}
\begin{split}
\frac{\underline{c}}{1 - \tgencc^{(p)}} \leq \tmix^{(p)} \leq \ceil{\tmix/p} \leq \ceil{\overline{c}/{p(1 - \tgencc)}}.
\end{split}
\end{equation*}
If $p \geq \tmix$, immediately, $1 - \tgencc^{(p)} \geq \underline{c} = 1/2$.
If $p \leq \tmix$, $ \ceil{\tmix/p} \leq 2 \frac{\tmix}{p}$, and the second claim follows. \qed

\subsection{Proofs of Theorem~\ref{theorem:fully-empirical-transition-kernel-learning}, Corollary~\ref{corollary:estimate-kappa}, and Lemma~\ref{lemma:asymptotic-behavior-matrix}}
\label{section:proof-empirical-estimation-kernel-infinity}
As part of our analysis, we will rely on the following result from the distribution learning literature, that
bounds with high probability the estimation error of learning a discrete distribution w.r.t total
variation from an iid sample, in terms of some measure of entropy of the empirical distribution.
\begin{theorem}[{\citet[Theorem~2.1]{cohen2020learning}}]
\label{theorem:empirical-distribution-learning}
Let $\delta \in (0,1), n \in \N, \mu \in \calP(\calX), X \sim \mu^{\otimes n}$, 
and 
$$\widehat{\mu} \eqdef \frac{1}{n} \sum_{t = 1}^{n}\sum_{x \in \calX} \pred{X_t = x} e_x,$$ 
the empirical measure. With probability at least $1 - \delta$, it holds that
$$\tv{\widehat{\mu} - \mu} \leq \frac{1}{\sqrt{n}}\nrm{\widehat{\mu}}_{1/2}^{1/2} + 3 \sqrt{\frac{\ln {2 /\delta}}{2n}}.$$
\end{theorem}

\begin{remark}
We note but do not pursue the fact that following \citet[Corollary~A.1]{cohen2020learning}, a more delicate analysis of the Rademacher estimates (involving the Wallis product for $\pi = 3.1415926\dots$) leads to a stronger constant for the first-order term,
$$\tv{\widehat{\mu} - \mu} \leq \sqrt{\frac{2}{\pi n}}\nrm{\widehat{\mu}}_{1/2}^{1/2} + \sqrt{\frac{2}{ \pi n^3}}\nrm{\widehat{\mu}^+}^{-1/2}_{-1/2} + 3 \sqrt{\frac{\ln {2 /\delta}}{2n}},$$
where 
$$\nrm{\widehat{\mu}^+}^{-1/2}_{-1/2} = \sum_{\substack{x \in \calX \\ \widehat{\mu}(x) > 0}}  1/\sqrt{\widehat{\mu}(x)}.$$
\end{remark}
From the law of total probability, for any $x \in \calX$,
\begin{equation}
\label{equation:decomposition-total-probability}
\begin{split}
\PR{\nrm{e_x(\widehat{P}_{\lambda} - P)}_1 > \widehat{W}_{\lambda, \delta}(X^m)} = \PR{\mathcal{E}_{0}} + \sum_{n_x = 1}^{m} \PR{\mathcal{E}_{n_x}},
\end{split}
\end{equation}
where for $n_x \in \set{0,1, \dots, m}$, we wrote
\begin{equation}
\label{equation:definition-of-errors}
\begin{split}
\mathcal{E}_{n_x} &\eqdef \set{ \nrm{e_x(\widehat{P}_{\lambda} - P)}_1 > \widehat{W}_{\lambda, \delta}(X^m) \text{ and } N_x = n_x}. \\
\end{split}
\end{equation}
We now treat the two cases where $n_x = 0$ and $n_x \neq 0$ separately.

\paragraph*{Case $n_x = 0$.}

In the event where $N_x = 0$, it also must be that $N_{x x'} = 0$ for any $x' \in \calX$, 
thus $\widehat{W}_{\lambda, \delta}(X^m) = 2$. 
Since $\nrm{e_x(\widehat{P}_{\lambda} - P)}_1 \leq 2$ almost surely, $\PR{\mathcal{E}_0} = 0$.

\paragraph*{Case $n_x \neq 0$.}
In this event, write $\widehat{P} = \widehat{P}_{\lambda}$ for $\lambda = 0$, 
the tally matrix without smoothing.
We can verify that
\begin{equation*}
\begin{split}
\nrm{e_x(\widehat{P}_{\lambda} - P)}_1 \leq \frac{N_x}{N_x +  \lambda(x) \abs{\calX}}
\nrm{e_x(\widehat{P} - P)}_1 + \frac{2 \lambda(x) \abs{\calX} }{N_x + \lambda(x)\abs{\calX}}.
\end{split}
\end{equation*}
Therefore, by algebraic manipulation and inclusion of events,
\begin{equation*}
\begin{split}
\PR{\mathcal{E}_{n_x}} \leq \PR{ \frac{1}{2}\nrm{e_x(\widehat{P} - P)}_1 > \frac{\nrm{e_x \widehat{P}}_{1/2}^{1/2}}{\sqrt{N_x}} + 3 \sqrt{\frac{\ln {(2 m \abs{\calX} /\delta)}}{2 N_x}} \text{ and } N_x = n_x}. \\
\end{split}
\end{equation*}
Similarly to \citet[p.9]{wolfer2021}, we use the technique credited to \citet[p.19]{billingsley1961statistical},
of simulating a trajectory from the Markov chain using an infinite array of independent random variables.
This allows us to write
\begin{equation*}
\begin{split}
\PR{\mathcal{E}_{n_x}} \leq \PR{\tv{e_x(\widetilde{P} - P)} > \frac{\nrm{e_x \widetilde{P}}_{1/2}^{1/2}}{\sqrt{n_x}} + 3 \sqrt{\frac{\ln {(2 m \abs{\calX} /\delta)}}{2 n_x}}} \leq \frac{\delta}{m \abs{\calX}}, \\
\end{split}
\end{equation*}
where $e_x\widetilde{P}$ corresponds to the empirical measure obtained by drawing independently
$$\widetilde{X}_1, \widetilde{X}_2, \dots, \widetilde{X}_{n_x} \sim e_x P^{\otimes n_x},$$
and the second inequality is by Theorem~\ref{theorem:empirical-distribution-learning}.
Plugging in \eqref{equation:decomposition-total-probability},
and applying two union bounds finishes proving the theorem.
\qed

Finally, to prove Lemma~\ref{lemma:asymptotic-behavior-matrix}, note that
by the ergodic theorem (Strong Law of Large Numbers, e.g. \citet[Theorem~C.1]{levin2009markov}), for any $x,x' \in \calX$,
\begin{equation*}
\begin{split}
N_x/m \as \pi(x), \qquad N_{x x'}/m \as  \pi(x) P(x,x').
\end{split}
\end{equation*}
The claims for $\widehat{W}_{\lambda, \delta}(X^m)$ then hold from 
algebraic manipulation,
the continuous mapping theorem,
and as the maximum over a finite collection of continuous functions is continuous.
\qed

\subsection{Proofs of Lemma~\ref{lemma:explicit-empirical-intervals} and Theorem~\ref{theorem:empirical-estimation-tmix}}
The proof of statement $(a)$ of Lemma~\ref{lemma:explicit-empirical-intervals} mirrors the steps of \citet[Theorem~2]{pmlr-v117-wolfer20a}, and is repeated below for the convenience of the reader. The distinction lies in relying on the sparsity-sensitive confidence interval of Theorem~\ref{theorem:fully-empirical-transition-kernel-learning}.
Successively,
\begin{equation*}
\begin{split}
&\PR{\abs{\esttgencc - \tgencc} > \frac{1}{S} + \max_{r \in [S]} \set{\frac{1}{r} \widehat{W}_{\lambda, \delta/S}^{(r)}(X^m)}} \\
&\stackrel{(i)}{\leq} \PR{\max_{s \in [S]}  \abs{\frac{\estcc_s - \cc_s}{s} }  > \max_{r \in [S]} \set{ \frac{1}{r} \widehat{W}_{\lambda, \delta/S}^{(r)}(X^m)}} \\
&\stackrel{(ii)}{\leq} \sum_{s=1}^{S}\PR{\abs{\estcc_s - \cc_s }  > s \max_{r \in [S]} \set{ \frac{1}{r} \widehat{W}_{\lambda, \delta/S}^{(r)}(X^m) }}   \\
&\stackrel{(iii)}{\leq} \sum_{s=1}^{S} \PR{ \nrm{\widehat{P}^{(s)}_{\lambda} - P^s}_\infty > \widehat{W}_{\lambda, \delta/S}^{(s)}(X^m) }   \\
&\stackrel{(iv)}{\leq} \sum_{s=1}^{S} \frac{\delta}{S} = \delta,
\end{split}
\end{equation*}
where $(i)$ follows from the fact that $\abs{\tgencc - \tgenccS} \leq 1/S$ 
and that for $\nu, \theta \in \R^S $ it is the case from sub-additivity of the uniform norm 
that $\abs{\nrm{\nu}_\infty - \nrm{\theta}_\infty} \leq \nrm{\nu - \theta}_\infty$;
$(ii)$ is an application of the union bound, $(iii)$ stems from the fact that 
the $\ell_\infty$ operator norm dominates the distance 
between Dobrushin contraction coefficients \citep[Fact~5.1]{pmlr-v117-wolfer20a}, 
and $(iv)$ is Theorem~\ref{theorem:fully-empirical-transition-kernel-learning}. 
We now show that modulo a logarithmic cost in the confidence intervals, 
we are allowed to randomize the subset of integers to explore and 
consider $\widehat{S}(X^m)$ as defined in \eqref{eq:adaptive-prefix} instead of a fixed $S \in \N$.
Recall from \eqref{eq:adaptive-prefix} that,
\begin{equation*}
\begin{split}
\widehat{S}(X^m) &= \ceil*{\frac{1}{ \widehat{K}^\phi_{\lambda, \delta}(X^m)}},
\end{split}
\end{equation*}
hence $\widehat{S}(X^m)$, by definition, take an integer value.
Furthermore,
\begin{equation*}
\begin{split}
    \widehat{K}^\phi_{\lambda, \delta}(X^m) &= \max_{s \in [\phi(m)]} \set{ \frac{1}{s} \widehat{W}^{(s)}_{\lambda, \delta}(X^m) } \geq \widehat{W}_{\lambda, \delta}(X^m).
\end{split}
\end{equation*}
We proceed to lower bound the right hand side. Assuming $m \geq 1$, at least one state $x_0 \in \calX$ was visited,
\begin{equation*}
\begin{split}
\widehat{W}_{\lambda, \delta}(X^m) &= 2 \max_{x \in \calX} \set{ \frac{  \sum_{x' \in \calX} \sqrt{N_{xx'}} + (3/\sqrt{2}) \sqrt{N_x} \sqrt{\ln {(2 m \abs{\calX} /\delta)}}  + \lambda(x)\abs{\calX} }{N_x + \lambda(x)\abs{\calX} }} \\
&\geq \frac{2}{\sqrt{m}}  \frac{\overbrace{\sum_{x' \in \calX} \sqrt{N_{x_0 x'}/N_{x_0}}}^{\geq 0} + \overbrace{(3/\sqrt{2}) \sqrt{\ln {(2 m \abs{\calX} /\delta)}}}^{\geq 1}  + \lambda(x_0)\abs{\calX}/\sqrt{N_{x_0}} }{\underbrace{\sqrt{N_{x_0}/m}}_{\leq 1} + \lambda(x_0)\abs{\calX}/\sqrt{m N_{x_0}} } \\
&\geq \frac{2}{\sqrt{m}}  \frac{ 1  + \lambda(x_0)\abs{\calX}/\sqrt{N_{x_0}} }{1 + \lambda(x_0)\abs{\calX}/\sqrt{m N_{x_0}} } \geq \frac{1}{\sqrt{m}},
\end{split}
\end{equation*}
thus almost surely $\widehat{S}(X^m) \in \set{1, \dots, \ceil{\sqrt{m}}}$.
As a consequence, 
\begin{equation*}
\begin{split}
&\PR{\abs{\esttgenccexplest - \tgencc} > \frac{1}{\widehat{S}} + \max_{r \in [\widehat{S}]} \set{\frac{1}{r} \widehat{W}_{\lambda, \delta/(\ceil{\sqrt{m}}\widehat{S} )}^{(r)}(X^m)}} \\
&\stackrel{(i)}{=}\sum_{S = 1}^{\ceil{\sqrt{m}}} \PR{\abs{\esttgenccexplest - \tgencc} > \frac{1}{\widehat{S}} + \max_{r \in [\widehat{S}]} \set{\frac{1}{r} \widehat{W}_{\lambda, \delta/(\ceil{\sqrt{m}}\widehat{S} )}^{(r)}(X^m)} \text{ and } \widehat{S}(X^m) = S} \\
&\stackrel{(ii)}{\leq}\sum_{S = 1}^{\ceil{\sqrt{m}}} \PR{\abs{\esttgenccexplestbb - \tgencc} > \frac{1}{S} + \max_{r \in [S]} \set{\frac{1}{r} \widehat{W}_{\lambda, \delta/(\ceil{\sqrt{m}} S )}^{(r)}(X^m)}} \\
&\stackrel{(iii)}{\leq}\sum_{S = 1}^{\ceil{\sqrt{m}}} \frac{\delta}{\ceil{\sqrt{m}}} = \delta \\
\end{split}
\end{equation*}
where $(i)$ stems from the law of alternatives and from the fact that 
$\widehat{S}(X^m) \in \set{1, \dots, \ceil{\sqrt{m}}}$ almost surely, 
$(ii)$ follows from inclusion of events,
and $(iii)$ is Lemma~\ref{lemma:explicit-empirical-intervals}-$(a)$ for a fixed $S \in \N$.
Invoking Corollary~\ref{corollary:control-tmix-with-kgen-alt2020}, 
$\tmix$ and $1/(1 - \tgencc)$ are equal up to the universal multiplicative constants $\underline{c}, \overline{c}$,
and plugging-in directly yields Theorem~\ref{theorem:empirical-estimation-tmix}.
\qed

\subsection{Proof of Theorem~\ref{theorem:learning-s-step-transition-matrix}}
We first treat the case $s=1$. 
By chaining of events,
\begin{equation*}
\begin{split}
\PR{\nrm{\widehat{P}_{\lambda} - P}_{\infty} > \eps} \leq \PR{\nrm{\widehat{P}_{\lambda} - P}_{\infty} > \widehat{W}_{\lambda, \delta/2}(X^m)} + \PR{\widehat{W}_{\lambda, \delta/2}(X^m) > \eps},
\end{split}
\end{equation*}
with $\widehat{W}_{\lambda, \delta/2}$ as defined in \eqref{equation:w-confidence-interval}.
The first term is smaller than $\delta/2$ as a consequence of Theorem~\ref{theorem:fully-empirical-transition-kernel-learning},
thus after an application of the union bound for $\widehat{W}_{\lambda, \delta/2}(X^m)$, it remains to control for $x \in \calX$,
\begin{equation*}
\begin{split}
\PR{ \frac{  \sum_{x' \in \calX} \sqrt{N_{xx'}} + (3/\sqrt{2}) \sqrt{N_x} \sqrt{\ln {(4 m \abs{\calX} /\delta)}}  + \lambda(x)\abs{\calX} }{N_x + \lambda(x)\abs{\calX} } > \eps/2}.
\end{split}
\end{equation*}
By the law of alternatives, and the simulation method 
of \citet[p.19]{billingsley1961statistical} on the event where $N_x  = n \in [m]$, 
we can reduce the above to the deviation of a function of independent 
random variables, $\widetilde{X}_1, \widetilde{X}_2, \dots, \widetilde{X}_n \sim e_x P^{\otimes n}$
\begin{equation}
\label{eq:deviation-independent}
\begin{split}
\PR{ \frac{  \sum_{x' \in \calX} \sqrt{\sum_{t=1}^{n} \pred{\widetilde{X}_t = x'}} + (3/\sqrt{2}) \sqrt{n} \sqrt{\ln {(4 m \abs{\calX} /\delta)}}  + \lambda(x)\abs{\calX} }{n + \lambda(x)\abs{\calX} } > \eps/2}.
\end{split}
\end{equation}
We first consider the case $n \in [m \pi(x)/2, 3 m \pi(x)/2]$, i.e. when $n$ is of the order of $\E[\pi]{N_x}$.
Writing 
$$\Psi(\widetilde{X}_1, \dots, \widetilde{X}_n) \eqdef \frac{  \sum_{x' \in \calX} \sqrt{\sum_{t=1}^{n} \pred{\widetilde{X}_t = x'}}}{n + \lambda(x)\abs{\calX} },$$
we bound its expectation,
\begin{equation}
\label{eq:expectation-by-line}
\begin{split}
\E{\Psi(\widetilde{X}_1, \dots, \widetilde{X}_n)} &\stackrel{(i)}{\leq}  \frac{1}{n + \lambda(x)\abs{\calX}} \sum_{x' \in \calX}\E{\sqrt{\Bin(n, P(x,x'))}} \\
&\stackrel{(ii)}{\leq}  \frac{\sqrt{n}}{n + \lambda(x)\abs{\calX}}\sum_{x' \in \calX}\sqrt{  P(x,x')} \leq  \frac{1}{\sqrt{n}} \nrm{e_x P}_{1/2}^{1/2} \\
\end{split}
\end{equation}
where $(i)$ is by linearity of the expectation, $(ii)$ follows from Jensen's 
inequality and the expression of the mean of a binomial random variable.
When $n \geq \frac{16 \nrm{e_x P}_{1/2}}{\eps^2}$, 
the expectation at \eqref{eq:expectation-by-line} is bounded by $\eps/4$.
For this regime, it is therefore enough that $m \geq \frac{32}{\eps^2}  \frac{\nrm{e_x P}_{1/2}}{\pi(x)}$.
In the rest of the proof $c > 0$ as customary in analysis, will denote a universal constant, 
that is allowed to change at each occurrence.
For $m \geq \frac{c \lambda(x) \abs{\calX}}{\eps \pi(x)}$, 
and $m \geq \frac{c}{\eps^2 \pi(x)} \ln \frac{1}{\delta \eps \pi(x)}$,
\begin{equation*}
\begin{split}
\frac{(3/\sqrt{2}) \sqrt{n} \sqrt{\ln {(4 m \abs{\calX} /\delta)}}  + \lambda(x)\abs{\calX} }{n + \lambda(x)\abs{\calX} } \leq \frac{\eps}{8},
\end{split}
\end{equation*}
and we are left with bounding 
\begin{equation*}
\begin{split}
\sum_{n \in [m \pi(x)/2, 3 m \pi(x)/2]} \PR{ \abs{\Psi - \E{\Psi}} > \eps/8} 
&\stackrel{(i)}{\leq} \sum_{n \in [m \pi(x)/2, 3 m \pi(x)/2]} \exp \left(- c n \eps^2\right) \\
&\stackrel{(ii)}{\leq} \frac{c}{\eps^2} \exp \left(- (c m \pi(x) \eps^2)/2\right), \\
\end{split}
\end{equation*}
where $(i)$ is McDiarmid's inequality, 
after having verified that $\Psi$ is $(2/n)$-Lipschitz with respect to the Hamming distance,
and $(ii)$ stems from $t > 0 \implies t e^{-t} \leq e^{-t/2}$.
The latter is smaller than $\delta/(4\abs{\calX})$ for $m \geq \frac{c}{\eps^2 \pi(x)} \ln \frac{\abs{\calX}}{\delta \eps}$.
It remains to handle the event where $N_x \not \in [m \pi(x)/2, 3 m \pi(x) / 2]$.
To do so, we invoke
the following Chernoff-Hoeffding tail 
inequalities.
\begin{corollary}[adaptation of {\citet[Theorem~3.1]{chung2012chernoff}} with  Corollary~\ref{corollary:control-tmix-with-kgen-alt2020}]
\label{corollary:chernoff-hoeffding-markov}
There exists universal constants $c_G, c_H$ such that the following holds.
For any $X = (X_1, \dots, X_m) \sim (\mu, P)$ with generalized contraction coefficient $\kappa_\star$ and stationary distribution $\pi$, and for any $\eta \in (0,1)$ and $f \colon \calX \to [0,1]$,
\begin{equation*}
\begin{split}
    \PR[\mu]{\frac{1}{m}\sum_{t = 1}^{m}f(X_t) > (1 + \eta)\bbE_{\pi} f(X) } &\leq c_G \nrm{\mu/\pi}_{\pi}\exp \left( -  \frac{1}{c_H}m \bbE_{\pi} f(X)  (1 - \kappa_\star) \eta^2  \right), \\
    \PR[\mu]{\frac{1}{m}\sum_{t = 1}^{m}f(X_t) < (1 - \eta)\bbE_{\pi} f(X) } &\leq c_G \nrm{\mu/\pi}_{\pi} \exp \left( -  \frac{1}{c_H} m \bbE_{\pi} f(X)  (1 - \kappa_\star) \eta^2  \right). \\
\end{split}
\end{equation*}
\end{corollary}
We readily obtain
\begin{equation}
\label{eq:concentration-simple}
\begin{split}
\PR[\mu]{N_x < \frac{m\pi(x)}{2}  } \leq \frac{c}{\pimin} \exp \left(- c m \pi(x)(1 - \kappa_\star) \right).
\end{split}
\end{equation}
We can handle general $s \in [m-1]$ similarly,
\begin{equation*}
\begin{split}
\PR[\mu]{N_x^{(s)} \leq \frac{1}{2} \ceil{(m-1)/s} \pi(x)} &\stackrel{(i)}{\leq} \frac{c}{\pimin} \exp \left(- c' \floor{(m-1)/s} \pi(x) \left(1 - \kappa_\star^{(s)}\right)\right) \\
&\stackrel{(ii)}{\leq} \frac{c}{\pimin} \exp \left(- c' (m - s - 1) \pi(x) (1 - \kappa_\star)\right) \\
\end{split}
\end{equation*}
where $(i)$ follows from $\E[\pi]{\pred{X_{1 + s(t-1)} = x}} = \pi(x)$, 
and $(ii)$ stems from 
Proposition~\ref{proposition:kgen-of-power}-$(a)$.
As a consequence for $m \geq c \frac{1}{\pi(x)(1 - \kappa_\star)} \ln \frac{1 }{\delta \pimin \eps}$, 
this error probability is smaller than $\frac{\delta}{2 \abs{\calX} \ceil{2/\eps}}$.
Taking a maximum over $x \in \calX$ and union bounds yield the theorem.
\qed

\subsection{Proof of Theorem~\ref{theorem:learn-tgencc-ub-absolute}}
\label{section:proof-point-estimator-absolute-error}
Following the same first steps $(i), (ii), (iii)$ as in the proof of Lemma~\ref{lemma:explicit-empirical-intervals} 
together with $S > \frac{2}{\eps}$,
\begin{equation*}
\begin{split}
\PR[\pi]{\abs{\esttgenccexpl - \tgencc} > \eps}
\leq \sum_{s=1}^{\ceil{2 /\eps}}\PR[\pi]{\nrm{\widehat{P}^{(s)} - P^{s}}_\infty  > s\frac{\eps}{2}}. \\
\end{split}
\end{equation*}
An application of Theorem~\ref{theorem:learning-s-step-transition-matrix} 
concludes the claim.
\qed

\subsection{Proof of Theorem~\ref{theorem:learn-tgencc-ub-relative}}
\label{section:proof-point-estimator-relative-error}
Let some free parameter $\eps \leq 1/8$ to be fixed later, and consider the estimator 
\begin{equation}
\label{eq:amplified-estimator}
\widehat{\kappa}^\times_\star \eqdef 1 - \frac{1 - \widehat{\kappa}^{(S)}}{\sqrt{c_\star} S},
\end{equation}
for some random integer $S$ defined by
\begin{equation*}
    \log_2 S \eqdef \argmin_{p \in \set{0,1, 2, \dots}} \set{ 1 - \widehat{\kappa}_\star^{(2^p)} > \frac{1}{4} + \eps},
\end{equation*}
where for $s \in \N$, $\widehat{\kappa}_\star^{(s)}$ is the estimator for $\kappa_s$ defined in Theorem~\ref{theorem:learn-tgencc-ub-absolute} that achieves additive precision $\eps$ with the skipped chain $X_{[s]}$ as input.
We define the events
\begin{equation*}
\begin{split}
    g(s, \eps) \eqdef \set{ \abs{ \widehat{\kappa}_\star^{(s)} - \kappa_\star^{(s)} } \leq \eps }, \qquad G(\eps) \eqdef \bigcap_{p = 0}^{\overline{p}} g(2^p, \eps),
\end{split}
\end{equation*}
where $\overline{p} \eqdef \ceil{\log_2 \tmix}$.
We consider a trajectory length of at least, $\max_{p \in \set{0, \dots, \overline{p}}} m_p$ with
\begin{equation*}
m_p \geq  c 2^p \max \set{ \frac{1}{\eps^2}\max_{s \leq \ceil{2/\eps}} \set{ \frac{\Gamma((P^{2^p})^s)}{s} }, \frac{1}{\eps^2 \pimin} \ln \frac{1}{\underline{\delta} \eps \pimin}, \frac{1}{\left(1 - \kappa_\star^{(2^p)}\right)\pimin}\ln \frac{1}{\pimin \underline{\delta}} },
\end{equation*}
with $\underline{\delta} \eqdef \frac{\delta}{1 + \overline{p}}$, i.e. from Proposition~\ref{proposition:kgen-of-power}, it suffices to take
\begin{equation*}
\begin{split}
m &\geq c \frac{1}{\eps^2} \max \set{ \Xi_\eps(P), \frac{1}{(1 - \kappa_\star)\pimin} \ln \frac{\log (1 - \kappa_\star)^{-1}}{\delta \eps \pimin}},
\end{split}
\end{equation*}
where 
\begin{equation*}
    \Xi_\eps(P) \eqdef \max_{p \in \set{0,1, \dots, \overline{p}}} \set{ 2^p \max_{s \leq \ceil{2 /\eps}} \set{ \frac{\Gamma\left(P^{s 2^p}\right)}{s}}}.
\end{equation*}
By the union bound and Theorem~\ref{theorem:learn-tgencc-ub-absolute},
\begin{equation*}
    \begin{split}
       \PR{G(\eps)^\complement} &\leq \sum_{p = 0}^{\overline{p}} \PR{g(2^p, \eps)^\complement} \leq (\overline{p} + 1) \underline{\delta} = \delta.
    \end{split}
\end{equation*}

\begin{lemma}
\label{lemma:good-event}
On the event $G(\eps)$,
\begin{enumerate}[label=$(\alph*)$]
    \item 
    The algorithm terminates before $S = 2^{\overline{p}}$.
    \item $1/4 < 1 - \kappa_\star^{(S)}$.
\end{enumerate}
\end{lemma}
\begin{proof}
Note that,
\begin{equation*}
\begin{split}
    1 - \kappa^{(2^{\overline{p}})}_\star &\geq 1 - \kappa_{\star}^{(\tmix)} \geq \frac{\underline{c}}{\tmix^{(\tmix)}} \geq \frac{\underline{c}}{\floor{\tmix / \tmix}} = \underline{c} = \frac{1}{2},
\end{split}
\end{equation*}
thus on $G(\eps)$,
\begin{equation*}
    1 - \widehat{\kappa}^{(2^{\overline{p}})}_\star > \frac{1}{2} - \eps > \frac{1}{4} + \eps,
\end{equation*}
and $(a)$ holds.
By the stopping condition of the algorithm, $1 - \widehat{\kappa}_\star^{(S)} > 1/4 + \eps$. Being on $G(\eps)$ implies $(b)$.
\end{proof}

Let
\begin{equation}
\label{eq:alpha-lower-bounds}
    \alpha > \max \set{ \frac{4 \eps + 1}{\sqrt{c_\star}} - 1
    , \frac{2 \overline{c} \sqrt{c_\star}}{1/4 + \eps} - 1 }.
\end{equation}
We proceed to verify the following inclusion of events,
\begin{equation*}
\begin{split}
    \set{\frac{1 - \widehat{\kappa}^\times_\star}{1 - \kappa_\star} > 1 + \alpha } \cap G(\eps) &\stackrel{(i)}{=} \set{\frac{1 - \widehat{\kappa}^{(S)}_\star}{\sqrt{c_\star} S(1 - \kappa_\star)} > 1 + \alpha  } \cap G(\eps) \\ 
    &\stackrel{(ii)}{\subset} \set{1 - \widehat{\kappa}^{(S)}_\star> (1 + \alpha) \sqrt{c_\star}\left(1 - \kappa_\star^{(S)}\right) } \cap G(\eps) \\
    &\subset \set{\abs{\kappa^{(S)}_\star - \widehat{\kappa}^{(S)}_\star}> \left(1 - \kappa^{(S)}_\star \right) \left[(1 + \alpha) \sqrt{c_\star} - 1 \right]  } \cap G(\eps) \\
    &\stackrel{(iii)}{\subset} \set{\abs{\kappa^{(S)}_\star - \widehat{\kappa}^{(S)}_\star}> \frac{1}{4} \left[(1 + \alpha) \sqrt{c_\star} - 1 \right]  } \cap G(\eps) \\
    &\subset \set{\abs{\kappa^{(S)}_\star - \widehat{\kappa}^{(S)}_\star}> \eps} \cap G(\eps) \\
    &\stackrel{(iv)}{\subset} \emptyset,
\end{split}
\end{equation*}
where $(i)$ is by definition of the estimator at \eqref{eq:amplified-estimator}, $(ii)$ follows from Proposition~\ref{proposition:kgen-of-power}-$(i)$,
$(iii)$ and $(iv)$ stem from Lemma~\ref{lemma:good-event}.
Furthermore, on $G(\eps)$, since $2^{\ceil{\log_2 x}} \leq 2x$, it holds that $S \leq 2 \tmix$. Thus,
\begin{equation*}
\begin{split}
    \set{\frac{1 - \kappa_\star}{1 - \widehat{\kappa}^\times_\star} > 1 + \alpha } \cap G(\eps)
    &\subset \set{\frac{2 \tmix (1 - \kappa_\star) \sqrt{c_\star}}{1 - \widehat{\kappa}^{(S)}_\star} > 1 + \alpha } \cap G(\eps)   \\
    &\stackrel{(i)}{\subset} \set{\frac{2 \overline{c} \sqrt{c_\star}}{1 - \widehat{\kappa}^{(S)}_\star} > 1 + \alpha  } \cap G(\eps) \\
    &\stackrel{(ii)}{\subset} \set{2 \overline{c} \sqrt{c_\star} > (1 + \alpha)(1/4 + \eps)} \cap G(\eps) \\
    &\subset \emptyset,
\end{split}
\end{equation*}
where $(i)$ stems from Corollary~\ref{corollary:control-tmix-with-kgen-alt2020} and $(ii)$ is by the definition of the stopping $S$.
Finally, we set $\eps = (\sqrt{2} - 1)/4 \in (0, 1/8)$. Then, from  \eqref{eq:alpha-lower-bounds} it is enough that $\alpha \geq 2\sqrt{2\log 4e} - 1$.

\qed

\subsection{Proof of Lemma~\ref{lemma:kgen-special-family}}

Recall the family of chains parametrized by $\eta \in (0, 2/3)$,
$P_\eta = \begin{psmallmatrix} 1/3 & 0 & 2/3 \\ \eta/2 & 1 - \eta & \eta/2 \\ \eta/2 & \eta/2 & 1 - \eta \end{psmallmatrix}$.
The spectrum and corresponding eigenvectors of $P_\eta$ are given by
\begin{alignat*}{3}
 &\lambda_1 = 1,\qquad  &&\lambda_2= \frac{1}{6} (2 - 3 \eta),\qquad  &&\lambda_3 = \frac{1}{2}(2 - 3 \eta), \\
 &v_1= (1, 1, 1),\qquad  &&v_2 = \left(-\frac{4}{3 \eta} , 1 , 1\right),\qquad  && v_3= \left(\frac{4}{4 - 9 \eta}, \frac{8 - 9 \eta}{9 \eta - 4}, 1 \right).
\end{alignat*}
We write $P_\eta = B_\eta D_\eta B_\eta \inv $
where
\begin{equation*}
\begin{split}
D_\eta = \begin{pmatrix} 1 & 0 & 0 \\ 0 & \frac{1}{6}(2 - 3 \eta) & 0 \\ 0 & 0 & \frac{1}{2}(2 - 3\eta) \end{pmatrix}, \qquad B_\eta = \begin{pmatrix} 1 & -\frac{4}{3 \eta} & \frac{4}{4 - 9 \eta}, \\ 1 & 1 & \frac{8 - 9 \eta}{9 \eta - 4} \\ 1 & 1 & 1 \end{pmatrix}, \\
B_\eta \inv = \begin{pmatrix} \frac{3 \eta}{3 \eta + 4} & \frac{4}{9 \eta + 12} & \frac{8}{9 \eta + 12}, \\ -\frac{3 \eta}{3 \eta + 4} & \frac{9 \eta^2}{2(9 \eta^2 + 6 \eta - 8)} & \frac{3 \eta(3 \eta - 4)}{2(9 \eta^2 + 6 \eta - 8)} \\ 0 & \frac{4 - 9 \eta}{6(3 \eta - 2)} & \frac{4 - 9 \eta}{12 - 18 \eta} \end{pmatrix}. \\
\end{split}
\end{equation*}
For $s \in \N$, $P_\eta^s = B_\eta D_\eta^s B_\eta \inv $
allows us to compute the
power chain,
as well as its contraction coefficient. 
\begin{equation*}
\begin{split}
\tv{e_1 P^s - e_2 P^s} &= \left( 1 - \frac{3 \eta}{2} \right)^{s -1} \left( 1 - \left( \frac{1 + 3^{-s}}{2} \right) \frac{3 \eta}{2} \right) \\
\tv{e_1 P^s - e_3 P^s} &= \left(1 - \frac{3 \eta}{2}\right)^{s-1} \frac{3^{-s}}{2} \left( \frac{3\eta}{2}  \left(\frac{3^s -1}{2}\right)   + 1 - \frac{3 \eta}{2} +  \left| \frac{3 \eta}{2} \left( \frac{3^s + 1}{2}\right) - 1\right| \right) \\
\tv{e_2 P^s - e_3 P^s} &= \left(1- \frac{3 \eta}{2} \right)^s \\
\end{split}
\end{equation*}
For $s \geq 1$, $1 + \frac{3 \eta}{2}(1 - 3^{-s}) \geq 0$, hence $\tv{e_1 P^s - e_2 P^s} \geq \tv{e_2 P^s - e_3 P^s}$. 
When $\eta > \frac{4}{3(3^s + 1)}$,
$$\tv{e_1 P^s - e_3 P^s} = \left(1 - \frac{3 \eta}{2} \right)^{s-1} \frac{3 \eta}{2} \left( \frac{1 - 3^{-s}}{2} \right) \leq \tv{e_1 P^s - e_2 P^s} .$$
When $\eta \leq \frac{4}{3(3^s + 1)}$, $\tv{e_1 P^s - e_3 P^s} = \left(1 - \frac{3 \eta}{2} \right)^{s} 3^{-s}$, that is also smaller than $\tv{e_1 P^s - e_2 P^s}$ from the considered ranges for $s$ an $\eta$. As a result,
$$(a) \qquad \cc_s(P_\eta) =\left( 1 - \frac{3 \eta}{2} \right)^{s -1} \left( 1 - \left( \frac{1 + 3^{-s}}{2} \right) \frac{3 \eta}{2} \right).$$
Claim $(b)$ follows directly from $(a)$ for $s = 1$;
for $s = 2$, $\cc_2 = (1 - 3 \eta / 2)(1 - 5 \eta / 6)$. Solving for $\eta$ yields $(c)$ from the assumption that $\eta < 1/6$.
\qed

\subsection{Proof of Lemma~\ref{lemma:special-category-chain-requires-arbitray-large-s}}

We compute the partial derivative of the function: 
$$\theta(s, \eta) \colon [1, \infty) \times (0, 1/6) \to \R,$$
\begin{equation*}
\begin{split}
\frac{\partial}{\partial s} \theta(s ,\eta) &= \frac{1-\left(1-\frac{3 \eta}{2}\right)^{s-1} \left(1-\frac{3}{4} \eta \left(3^{-s}+1\right)\right)}{s^2} \\
&+\frac{\frac{1}{4} (-1) \eta 3^{1-s} \log (3) \left(1-\frac{3 \eta}{2}\right)^{s-1}}{s} \\
&- \frac{\left(1-\frac{3}{4} \eta \left(3^{-s}+1\right)\right) \left(1-\frac{3 \eta}{2}\right)^{s-1} \log \left(1-\frac{3 \eta}{2}\right)}{s} \\
&= \frac{1}{s^2}\left( 1-\left(1-\frac{3 \eta}{2}\right)^{s-1}\left(1 -  h(s, \eta)\right) \right), \\
\end{split}
\end{equation*}
with 
\begin{equation}
\label{equation:critical-function}
\begin{split}
h(s, \eta) \eqdef  \frac{3}{4} \eta \left(3^{-s}+1\right) \left( 1 + s \log \left(\frac{1}{3}-\frac{\eta}{2}\right) \right) + \frac{3}{4} s \eta \log(3) - s  \log \left(1-\frac{3 \eta}{2}\right).
\end{split}
\end{equation}
Suppose that $1 \leq s \leq \ln \frac{1}{\eta}$. Then,
\begin{equation*}
\begin{split}
h(s, \eta) \geq \underline{h}(s, \eta) \eqdef \frac{3}{4} \eta \left(\eta^{\ln{3}}+1\right) \left( 1 + s \log \left(\frac{1}{3}-\frac{\eta}{2}\right) \right) + \frac{3}{4} \eta \log(3) - s  \log \left(1-\frac{3 \eta}{2}\right).
\end{split}
\end{equation*}
Taking the partial derivative, and as $x > 0 \implies 1 - 1/x \leq \ln x \leq x - 1$,
\begin{equation*}
\begin{split}
\frac{\partial }{\partial s} \underline{h}(s, \eta) &= \frac{3}{4} \eta \left(\eta^{\ln{3}}+1\right) \log \left(\frac{1}{3}-\frac{\eta}{2}\right) - \log \left(1-\frac{3 \eta}{2}\right) \\
&\geq \frac{3}{4} \eta \left(\eta + 1\right) \left( \log \frac{1}{3} + 1 - \frac{1}{1 - 3 \eta / 2} \right) + \frac{3 \eta}{2}, \\
\end{split}
\end{equation*}
which is positive for $\eta \in (0, 1/6)$.
It follows that $\underline{h}$ is an increasing function of $s$ for any fixed $\eta$. Taking the limit,
$$\lim_{s \rightarrow 1} \underline{h}(s, \eta) = \frac{\partial }{\partial s} \underline{h}(s, \eta) + \frac{3}{4} \eta \left(\eta^{\ln{3}}+1\right) + \frac{3}{4} \eta \log(3) \geq 0.$$
As a consequence, we successively have that for $\eta \in (0, 1/6)$ and $s \in [1, \ln (1 / \eta)]$, $h(s, \eta)$ is positive, $\frac{\partial}{\partial s} \theta(s ,\eta)$ is positive,
$s \mapsto \theta(s ,\eta)$ is increasing, and  
$$s_\star(\eta) \in \argmax_{s \in \N} \set{ \theta(s, \eta)} \implies s_\star(\eta) > \ln(1/\eta).$$ We illustrate this effect in Figure~\ref{fig:special-chain-eta}.
\qed

\begin{figure}
\begin{center}

\begin{tikzpicture}
\begin{axis}[
    legend cell align={left},
    width=12cm,
    height=3cm,
    axis lines = left,
    xlabel = $s$,
    ylabel = {$\frac{1}{s}\left( 1 - \kappa_s(P_\eta) \right)$},
]
\draw ({axis cs:4.968121337890633,0}|-{rel axis cs:0,0}) -- ({axis cs:4.968121337890633,0}|-{rel axis cs:0,1});
\addplot [
    domain=0:20, 
    samples=1000, 
    color=red,
    ]
    {(1 - ( 1 - 3 * (0.03) / 2 )^(x-1) * ( 1 - ((1 + 3^(-x))/2) * 3 * (0.03) / 2 )) / x};
\addlegendentry{$\eta = 0.03$}
\node[above] at (45,5) {$s_\star$};
\end{axis}
\end{tikzpicture}

\begin{tikzpicture}
\begin{axis}[
    legend cell align={left},
    width=12cm,
    height=3cm,
    axis lines = left,
    xlabel = $s$,
    ylabel = {$\frac{1}{s}\left( 1 - \kappa_s(P_\eta) \right)$}
]
\draw ({axis cs:6.065179347991954,0}|-{rel axis cs:0,0}) -- ({axis cs:6.065179347991954,0}|-{rel axis cs:0,1});
\addplot [
    domain=0:20, 
    samples=1000, 
    color=red,
    ]
    {(1 - ( 1 - 3 * (0.02) / 2 )^(x-1) * ( 1 - ((1 + 3^(-x))/2) * 3 * (0.02) / 2 )) / x};
\addlegendentry{$\eta = 0.02$}
\node[above] at (55,5) {$s_\star$};
\end{axis}
\end{tikzpicture}

\begin{tikzpicture}
\begin{axis}[
    legend cell align={left},
    width=12cm,
    height=3cm,
    axis lines = left,
    xlabel = $s$,
    ylabel = {$\frac{1}{s}\left( 1 - \kappa_s(P_\eta) \right)$},
]
\draw ({axis cs:8.48388725519182,0}|-{rel axis cs:0,0}) -- ({axis cs:8.48388725519182,0}|-{rel axis cs:0,1});
\addplot [
    domain=0:20, 
    samples=1000, 
    color=red,
    ]
    {(1 - ( 1 - 3 * (0.01) / 2 )^(x-1) * ( 1 - ((1 + 3^(-x))/2) * 3 * (0.01) / 2 )) / x};
\addlegendentry{$\eta = 0.01$}
\node[above] at (80,5) {$s_\star$};
\end{axis}
\end{tikzpicture}

\caption{Decreasing $\eta$ moves the mode -- represented by a vertical line -- to the right.}
\label{fig:special-chain-eta}
\end{center}
\end{figure}

\section*{Acknowledgments}
The author thanks anonymous referees for their insightful comments which helped improve the quality of this manuscript. 
Special thanks to an anonymous referee who suggested investigating amplification methods for this problem, and Pierre Alquier for the interesting conversations.

\bibliography{bibliography}
\bibliographystyle{abbrvnat}

\newpage

\appendix
\section{Numerical tightness of  Theorem~\ref{theorem:control-tau-with-kgen-and-xi-alt2020} and Theorem~\ref{theorem:control-tau-with-kgen-and-xi-converse}}
\label{section:numerical-tightness-bound}

To illustrate Theorem~\ref{theorem:control-tau-with-kgen-and-xi-alt2020} and Theorem~\ref{theorem:control-tau-with-kgen-and-xi-converse}, we proceed with the following simulation.
We define a probability distribution $\calD(\calW(\calX))$ on the set of transition matrices over $\calX$, 
where $P \sim \calD(\calW(\calX))$ means that for any $x \in \calX$, the $x$th row of the random transition
matrix is sampled from a Dirichlet distribution with uniform concentration parameters:
 $e_x P \sim \Dir( 1/ \abs{\calX})$.
In Figure~\ref{fig:proxy-universal-constant}, we take a random support size $\abs{\calX} \sim \Uniform \set{3 \dots 8}$, and then sample $n = 500$ pairs $(P_k, \xi_k)$
where $\xi_k \sim \Uniform([0, 1/2])$ and $P_k \sim \calD(\calW(\calX))$. 
For each pair, we plot $\tau(\xi_k)$ as a function of $\xi_k$, the upper and lower 
bounds at Theorem~\ref{theorem:control-tau-with-kgen-and-xi-alt2020}, 
as well as the function $\xi \to \frac{\ln{(e/(2\xi))}}{1 - \tgencc}$. 
Our experiments invite us to conjecture that Theorem~\ref{theorem:control-tau-with-kgen-and-xi-alt2020}-$(a)$ actually holds with the stronger constant $c = 1/2$.

\begin{figure}
\centering
\begin{tikzpicture}
\begin{axis}[
    legend cell align={left},
    width=12cm,
    axis lines = left,
    xlabel = $\xi$,
    ymin = 0, ymax=6
]

\addplot [
    domain=0:0.5, 
    samples=100, 
    color=red,
]
{1 - x};
\addlegendentry{$(1 - \xi)$}

\addplot [
    domain=0:0.5, 
    samples=100, 
    color=green,
    ]
    {ln(e / x)};
\addlegendentry{$\ln(c e / \xi)$ with $c = 1$}

\addplot [
    domain=0:0.5, 
    samples=100, 
    color=blue,
    ]
    {ln((1/2)* e / x)};
\addlegendentry{$\ln(c e / \xi)$ with $c = 1/2$}

\addplot[
    only marks,
    mark=*,
    mark options={black},
    mark size=1.0pt]
table[meta=xi]
{plot-generalized-contraction-coefficient.dat};
\addlegendentry{$\tau(\xi)(1 - \tgencc)$ for $P_1, P_2, \dots P_{500} \sim \calD(\calW(\calX))$}

\end{axis}
\end{tikzpicture}

\caption{Mixing time $\tau(\xi)$, with upper and lower bounds of 
Theorem~\ref{theorem:control-tau-with-kgen-and-xi-alt2020} for a random sample of transition kernels over $\calX, \abs{\calX} \sim \Uniform\set{3 \dots 8}$.}
\label{fig:proxy-universal-constant}
\end{figure}

\end{document}